\documentclass[titlepage,final,12pt]{article}
\usepackage{graphics}
\usepackage{graphicx}
\usepackage{subfigure}
\usepackage{epsfig}
\usepackage{enumerate}
\usepackage{amsfonts}
\usepackage{theorem}
\usepackage{array}
\usepackage[reqno]{amsmath}
\usepackage{color}

\theoremstyle{change}
\newtheorem{proclaim}{PROCLAIM}[section]
\newtheorem{theorem}[proclaim]{Theorem}

\newtheorem{proposition}[proclaim]{Proposition}
\newtheorem{corollary}[proclaim]{Corollary}

\def\state #1. { \noindent{\bf#1.\enspace}}

\newcommand{\comp}{\,{\raise 1pt \hbox{$\scriptstyle\circ$}}\,}

\newcommand{\con}{\mathop{\rm con}}

\newcommand{\dist}{\mathop{\rm dist}}
\newcommand{\exs}{\mathop{\rm exs}}
\newcommand{\epi}{\mathop{\rm epi}}

\newcommand{\dom}{\mathop{\rm dom}}

\newcommand{\lev}{\mathop{\rm lev}}

\newcommand{\reals}{\mathbb{R}}
\newcommand{\preals}{\mathbb{R}_+}
\newcommand{\Reals}{\overline{\mathbb{R}}}

\newcommand{\natnums}{{{\rm l} \kern -.13em {\rm N} }}
\newcommand{\nats}{\mathbb{N}}
\newcommand{\snats}{{I\kern -.29em N}}
\newcommand{\rats}{{Q\kern -.64em \raise 1pt \hbox{$\scriptstyle |$}\;\,}}
\newcommand{\srats}
    {{Q\kern -.56em \raise 1.2pt \hbox{$\scriptscriptstyle /$}\,}}
\newcommand{\ints}{Z\kern -.46em Z}
\newcommand{\ball}{\mathbb{B}}
\newcommand{\Ex}{{I\kern-.35em E}}
\newcommand{\pluss}{\hskip1pt \raise1pt\vbox{\hrule width6pt \vskip1pt \hrule
                    width6pt} \kern-4pt{\lower1pt\hbox{\vrule height6pt
            \kern1pt\vrule height6pt}}\hskip5pt}
\newcommand{\eop}
    {\hfill{$\vcenter{\hrule height1pt \hbox{\vrule width1pt height5pt
     \kern5pt \vrule width1pt} \hrule height1pt}$} \medskip}

\newcommand{\setd}{{ d \kern -.15em l}}
\newcommand{\hatsetd}{ d \hat{\kern -.15em l }}

\renewcommand{\epsilon}{\varepsilon}
\renewcommand{\phi}{\varphi}

\newcommand{\tto}{\;{\lower 1pt \hbox{$\rightarrow$}}\kern -12pt
           \hbox{\raise 2.5pt \hbox{$\rightarrow$}}\;}
\newcommand{\overto}[1]{\,{\raise 0pt\hbox{$\rightarrow$}}\kern -9pt
     \hbox{\lower 3pt \hbox{$\scriptscriptstyle#1$}}\hskip6pt}
\newcommand{\underto}[1]{\,{\lower 1pt\hbox{$\rightarrow$}}\kern -9pt
     \hbox{\raise 4pt \hbox{$\,\scriptscriptstyle#1$}}\hskip7pt}
\newcommand{\bigoverto}[1]{{\raise 0pt\hbox{$\,\longrightarrow$}}\kern -16pt
     \hbox{\lower 3pt \hbox{$\scriptscriptstyle#1$}}\hskip4pt}
\newcommand{\bigunderto}[1]{\,{\lower 1pt\hbox{$\longrightarrow$}}\kern -16pt
     \hbox{\raise 4pt \hbox{$\,\scriptscriptstyle#1$}}\hskip6pt}
\newcommand{\bigbigto}[2]{\,{\raise 0pt\hbox{$\,\longrightarrow$}}\kern -16pt
     \hbox{\lower 3pt \hbox{$\scriptscriptstyle#2$}}\kern -10pt
     \hbox{\raise 4pt \hbox{$\,\scriptscriptstyle#1$}}\hskip7pt}
\newcommand{\downto}{{\raise 1pt \hbox{$\scriptstyle \,\searrow\,$}}}
\newcommand{\upto}{{\raise 1pt \hbox{$\scriptstyle \,\nearrow\,$}}}

\newcommand{\notimply}
    {\quad\hbox{$\Longrightarrow \kern -14pt {/}$}\hskip6pt\quad}

\newcommand{\low}[1]{{\lower1pt \hbox{$\scriptstyle #1$}}}
\newcommand{\loww}[1]{{\lower2pt \hbox{$\scriptstyle #1$}}}
\newcommand{\high}[1]{{\raise1pt \hbox{$\scriptstyle #1$}}}

\newcommand{\nlev}{\mathop{\lev}\nolimits}

\newcommand{\ninf}{\mathop{\rm inf}\nolimits}
\newcommand{\nsup}{\mathop{\rm sup}\nolimits}
\newcommand{\nmax}{\mathop{\rm max}\nolimits}

\newcommand{\nnmin}{\mathop{\rm minimize}}

\newcommand{\nargmin}{\mathop{\rm argmin}\nolimits}

\newcommand{\lwdy}[2]{\mathrel{\mathop
        {\raisebox{0.1ex}{\null$#1$}}{\hbox{\kern -1.0em
    {\raisebox{-0.8ex}{$\scriptstyle{\;\to #2}$}}}}}}
\newcommand{\lwwdy}[2]{\mathrel{\mathop
        {\raisebox{0.2ex}{\null$#1$}}{\hbox{\kern -1.0em
    {\raisebox{-1.1ex}{$\scriptstyle{\;\to #2}$}}}}}}
\newcommand{\slwwdy}[2]{\scriptsize{{\mathrel{\mathop
        {\raisebox{0.2ex}{\null$#1$}}{\hbox{\kern -1.0em
    {\raisebox{-1.1ex}{$\scriptstyle{\;\to #2}$}}}}}}}}

\def\eto{\,{\lower 1pt\hbox{$\rightarrow$}}\kern -11pt
     \hbox{\raise 4pt \hbox{$\, \scriptstyle e$}}\hskip7pt}
\def\hto{\,{\lower 1pt\hbox{$\rightarrow$}}\kern -11pt
     \hbox{\raise 4pt \hbox{$\, \scriptstyle h$}}\hskip7pt}
\def\pto{\,{\lower 1pt\hbox{$\rightarrow$}}\kern -11pt
     \hbox{\raise 4.5pt \hbox{$\, \scriptstyle p$}}\hskip7pt}
\def\cto{\,{\lower 1pt\hbox{$\rightarrow$}}\kern -11pt
     \hbox{\raise 4pt \hbox{$\, \scriptstyle c$}}\hskip7pt}

\def\dy#1\until#2{\mathrel{\mathop
        {\null#1}\limits^{\hbox{\lower 1.3ex \hbox{$\scriptstyle{\;\to}$}}}}
        {\hbox{{\kern -.2em${\null}
       ^{\hbox{\raise .2ex \hbox{$\scriptstyle{#2}$}}}$}}}}
\def\hidy#1\until#2{\mathrel{\mathop
        {\null#1}\limits^{\hbox{\lower 1.0ex \hbox{$\scriptstyle{\;\to}$}}}}
        {\hbox{{\kern -.2em${\null}
	^{\hbox{\raise 0.9ex \hbox{$\scriptstyle{#2}$}}}$}}}\!}
\def\lody#1\until#2{\mathrel{\mathop
        {\null#1}\limits_{\hbox{\raise 0.5ex \hbox{$\scriptstyle{\to}$}}}}
        {\hbox{{\kern -.25em${\null}
       _{\hbox{\lower 0.7ex \hbox{$\scriptstyle{#2}$}}}$}}}}
\def\lowdy#1\until#2{\mathrel{\mathop
        {\null#1}\limits_{\hbox{\raise 1.0ex \hbox{$\scriptstyle{\to}$}}}}
        {\hbox{{\kern -.25em${\null}
       _{\hbox{\lower 0.9ex \hbox{$\scriptstyle{#2}$}}}$}}}}

\def\gph{\mathop{\rm gph}\nolimits}

\newcommand{\bcdot}{\,{\raise .2ex \hbox{$\cdot$}}\,}

\begin{document}


\begin{center}
\begin{large}
{\bf Stability and Error Analysis for\\ Optimization and Generalized Equations}
\smallskip
\end{large}
\vglue 0.7truecm
\begin{tabular}{c}
  \begin{large} {\sl Johannes O. Royset 
                                  } \end{large} \\
  Operations Research Department\\
  Naval Postgraduate School\\
  joroyset@nps.edu\\
\end{tabular}

\vskip 0.2truecm

\end{center}

\vskip 1.3truecm

\noindent {\bf Abstract}. \quad Stability and error analysis remain challenging for problems that lack regularity properties near solutions, are subject to large perturbations, and might be infinite dimensional. We consider nonconvex optimization and generalized equations defined on metric spaces and develop bounds on solution errors using the truncated Hausdorff distance applied to graphs and epigraphs of the underlying set-valued mappings and functions. In the process, we extend the calculus of such distances to cover compositions and other constructions that arise in nonconvex problems. The results are applied to constrained problems with feasible sets that might have empty interiors, solution of KKT systems, and optimality conditions for difference-of-convex functions and composite functions.\\

\vskip 0.5truecm

\halign{&\vtop{\parindent=0pt
   \hangindent2.5em\strut#\strut}\cr
                         \hfill\break
\hglue 1.35cm approximation theory, perturbation analysis, solution stability.
                         \cr\cr

{\bf Date}:\quad \ \today \cr}

\baselineskip=15pt

\section{Introduction}

Since the early days of convex analysis, epigraphs have been central to understanding functions in the context of minimization problems. Local properties of epigraphs can be used to define subgradients while global properties characterize convexity and lower semicontinuity. The distance between two epigraphs bounds the discrepancy between the corresponding minima and near-minimizers. Likewise, set-valued mappings can be fully represented by their graphs, with graphical convergence being key to understanding approximations of solutions of generalized equations defined by such mappings. These {\it set-based perspectives} lead to a unified approach to stability and error analysis for a wide range of variational problems. In this paper, we estimate the {\it truncated Hausdorff distance} between sets
and demonstrate that it provides insight about the stability of constraint systems and optimization problems even when the feasible sets have empty interiors.
Without assuming any local properties, we establish that the truncated Hausdorff distance bounds the discrepancy between near-solutions of two generalized equations when applied to the graphs of the underlying set-valued mappings. The result is illustrated in the context of optimality conditions for difference-of-convex functions, composite functions, and nonlinear programs. Throughout, we focus on nonconvex problems. Most of the results are established for general metric spaces and therefore apply broadly, including in areas such as nonparametric statistics, optimal control, function identification, and decision rule optimization.

Stability  and error analysis for optimization and, more generally, variational problems have been developed from several angles; see for example \cite{Laurent.72,Attouch.84,Polak.97,VaAn,BonnansShapiro.00} for comprehensive treatments. There is an extensive literature on {\it local} stability based on metric regularity and calmness
\cite{IoffeOutrata.08,Penot.10}, tilt-stability
\cite{EberhardWenczel.12,LewisZhang.13,DrusvyatskiyLewis.13}, full-stability
\cite{MordukhovichRockafellarSarabi.13}, and connections with iterative schemes
\cite{KlatteKrugerKummer.12}; see also the monographs
\cite{AubinEkeland.84,Mordukhovich.13a,Mordukhovich.13b}
and the surveys \cite{Pang.97,Aze.03}. This paper takes an alternative, {\it global} perspective that can be traced back to the late 60s and pioneering studies of the truncated Hausdorff distance between convex cones \cite{WalkupWets.67} and general convex sets \cite{Mosco.69}. The full potential of the approach emerges in \cite{AttouchWets.91,AttouchWets.93a,AttouchWets.93b}, which establish that the truncated Hausdorff distances between epigraphs furnish bounds on the corresponding discrepancies between minima and minimizers; see also \cite{AzePenot.90,AtLW91,BeerLucchetti.91,BeerLucchetti.92} for parallel developments and especially the monograph \cite{Beer.93} with its detailed treatment of topologies and metrics on spaces of closed sets. From the myriad of possibilities the Attouch-Wets distance \cite{AttouchWets.86} emerges as the theoretically most useful by virtue of being a metric on spaces of nonempty closed sets as well as other factors. Still, we concentrate on the truncated Hausdorff distance due to its more intuitive form and direct relationship to quantities of interest such as minima and minimizers. It anyhow furnishes accurate estimates of the Attouch-Wets distance \cite{VaAn,Royset.18}. This global perspective based on set distances provides foundations for computationally attractive approximations of functions \cite{RoysetWets.15b,Royset.18,Royset.19} and formulations of function identification problems \cite{RoysetWets.15b}, especially in nonparametric statistics \cite{RstW13:1dim,RoysetWets.19}.

The difficulty of estimating the truncated Hausdorff distance for actual problem instances remains a major hurdle for its practical use. Fundamental results and calculus rules are laid out in \cite{AzePenot.90b,AttouchWets.91}, but mostly for epigraphs in the convex case. Results on epi-multiplication and epi-sums are given in \cite{AttouchWets.91}. Inverse images of convex sets are well-behaved under sufficiently small perturbations. This fact enables the development of results for intersections of sets and sums of functions  in the convex case \cite{AzePenot.90b}. Since the Legendre-Fenchel transform is an isometry for lower semicontinuous proper convex functions under a closely related pseudo-metric defined in terms of the epi-regularized functions \cite{AttouchWets.86}, additional estimates of the truncated Hausdorff distance emerge via the dual operations under this transform \cite{AttouchWets.91}. In this paper, we switch the focus to nonconvex sets and functions and develop a series of results that support calculations of the truncated Hausdorff distance in practice.

Section 2 lays out the terminology and provides some motivating facts. Section 3 develops estimates for the truncated Hausdorff distance between arbitrary sets. Section 4 turns to specific results for epigraphs and applications in disjunctive programming, formulations with constraint softening, and penalty methods. Section 5 extends the methodology to set-valued mappings and demonstrates its usefulness for generalized equations such as those arising from optimality conditions. An appendix supplements with proofs.

\section{Distances and Applications}

For a point $x$ in a metric space $(X,d_X)$ and $C\subset X$, we denote by $\dist(x,C)$ the usual {\it point-to-set distance}, i.e.,
\begin{equation*}
\dist(x, C) := \inf\left\{ d_X(x,\bar x) ~|~ \bar x\in C \right\} \mbox{ if $C$ is nonempty and } \dist(x,\emptyset) := \infty.
\end{equation*}
The {\it excess} of $C$ over $D\subset X$ is given by
\begin{equation*}
  \exs(C;D) := \nsup\{ \dist(x, D)~|~x\in C\} \mbox{ if } C,D \mbox{ are nonempty},
\end{equation*}
$\exs(C;D) := \infty$ if $C$ nonempty and $D$ empty, and $\exs(C;D) := 0$ otherwise.
The Pompeiu-Hausdorff distance between $C$ and $D$ is $\max\{ \exs(C;D), \exs(D;C)\}$, but tends to be infinity for unbounded sets and therefore is not central to our development. Instead, we rely on a localization argument relative to a point $x^{\text{ctr}} \in X$, which we call the
{\it centroid} of $X$. The choice of centroid can be made
arbitrarily, but results might be sharper if it is near the ``interesting'' parts of the sets at hand as we often restrict the attention to intersections of sets with the centered closed ball
\[
\ball_X(\rho) := \{x\in X~|~d_X(x^{\text{ctr}}, x) \leq \rho\} \mbox{ for } \rho\geq 0.
\]
Given $\rho\geq 0$, we define the {\it  truncated Hausdorff distance} between two sets $C,D\subset X$ as
\begin{equation*}
  \hatsetd_\rho(C,D) := \max\Big\{\exs\big(C \cap \ball_X(\rho); D\big), ~\exs\big(D \cap \ball_X(\rho); C\big)\Big\},
\end{equation*}
which is always finite as long as $C$ and $D$ are nonempty and $\rho<\infty$.
Trivially, $\hatsetd_\infty(C,D)$ is the Pompeiu-Hausdorff distance between $C$ and $D$, but we focus on finite $\rho$ in the following.

The notation for the truncated Hausdorff distance suppresses its dependence on the choice of metric and centroid. The following results holds for all metrics and centroids unless otherwise specified. In particular,\\
\begin{quote}
for a normed linear space the metric is consistently assumed to be the one induced by the norm and the centroid is the zero point of the space.\\
\end{quote}
This is a harmless assumption, easily overcome, but kept here to simplify expressions. The ``hat-notation'' hints to a broader landscape of closely related distances between sets including the Attouch-Wets metric; see \cite[Chapter 4]{VaAn} for a summary of results. Although the truncated Hausdorff distance fails to be a metric on spaces of nonempty closed sets, it is obviously nonnegative and symmetric. A triangle inequality of sort also holds. Let $\preals := [0,\infty)$.

\begin{proposition}\label{prop:triangle}{\rm (triangle inequality, extended sense).}
  For a metric space $X$ with centroid $x^{\rm ctr}$, sets $C_1, C_2, C_3\subset X$, and $\rho\in \preals$,
  \[
  \hatsetd_\rho(C_1, C_3) \leq \hatsetd_{\bar\rho}(C_1,C_2) + \hatsetd_{\bar\rho}(C_2,C_3)
  \]
  provided that $\bar\rho> 2\rho + \max_{i=1, 2, 3} \dist(x^{\rm ctr}, C_i)$.
\end{proposition}
\state Proof. The arguments in the proofs of \cite[Prop. 1.2]{AttouchWets.91} and \cite[Prop. 3.1]{Royset.18} can easily be modified for the present assumptions.\eop

For a function $f:X\to \Reals := [-\infty,\infty]$, the characterizing set in the context of minimization problems is its {\it epigraph}
\[
\epi f := \big\{(x,\alpha) \in X\times \reals~|~ f(x) \leq \alpha\big\}.
\]
The truncated Hausdorff distance between epigraphs requires a metric and centroid for $X\times \reals$ and we consistently adopt\\

\begin{quote}
the product metric $((x,\alpha),(\bar x, \bar \alpha))\mapsto \max\{d_X(x,\bar x), |\alpha-\bar \alpha|\}$ and centroid $(x^{\text{ctr}},0)$, where $x^{\text{ctr}}$ is a centroid of $X$.\\
\end{quote}

The main motivation for studying the truncated Hausdorff distance between epigraphs is its relation to minima and minimizers. We recall that $\inf f := \inf \{f(x)~|$ $x\in X\}$, $\epsilon\mbox{-}\nargmin f := \{x \in \dom f~|~f(x)\leq \inf f + \epsilon\}$ for $\epsilon\geq 0$, with $\dom f$ $:=$ $\{x\in X~|~f(x) < \infty\}$, and $\nlev_\delta f := \{x\in X~|~f(x)\leq \delta\}$ for $\delta \in \Reals$. (We adopt the usual  arithmetic rules for extended real-valued numbers with an orientation towards minimization so that $\infty - \infty$ as well as $- \infty + \infty$ are set to $\infty$; see \cite[1.E]{VaAn}.)
The application in the context of minimization problems becomes clear from the following two propositions, which are essentially in \cite{AttouchWets.93a,Royset.18}. Still, due to minor adjustments in assumptions we provide proofs in the appendix.

\begin{proposition}\label{prop:solestimates}{\rm (approximation of infima and near-minimizers).}
For a metric space $X$, functions $f,g:X\to \Reals$, and $\epsilon, \rho \in \preals$,
\begin{align*}
| \inf f - \inf g| & \leq \hatsetd_\rho(\epi f, \epi g)\\
\exs\big(\epsilon\mbox{-}\nargmin g \cap \ball_X(\rho); ~\delta\mbox{-}\nargmin f\big) &\leq \hatsetd_\rho(\epi f, \epi g)
\end{align*}
provided that $\ninf f, \inf g \in [-\rho, \rho-\epsilon)$ and $\gamma\mbox{-}\nargmin f \cap \ball_X(\rho)$ as well as $\gamma\mbox{-}\nargmin g \cap \ball_X(\rho)$ are nonempty for all $\gamma>0$, with the second assertion also requiring $\delta > \epsilon + 2\hatsetd_\rho(\epi f, \epi g)$.
\end{proposition}

These bounds are sharp as discussed in \cite{Royset.18}. We note that $\delta$ cannot generally be equal to $\epsilon + 2\hatsetd_\rho(\epi f, \epi g)$. For example, suppose that $f(x) = x$ for $x>0$ and $f(x) = \infty$ otherwise; and $g(x) = x$ for $x\geq 0$ and $g(x) = \infty$ otherwise. Then, for $\rho\geq 0$, $\hatsetd_\rho(\epi f, \epi g)=0$, $\nargmin g = \{0\}$, $\nargmin f = \emptyset$, and $\exs(\nargmin g; \nargmin f) = \infty$. The role of $\rho$ emerges from the proposition: it needs to be large enough so that the epigraphs intersected with $\ball_{X\times \reals}(\rho)$ retain points corresponding to infima and near-minimizers.

\begin{proposition}\label{prop:levelsets}{\rm (approximation of level sets).} For a metric space $X$, functions $f,g:X\to \Reals$, $\rho\in \preals$, and $\delta \in [-\rho,\rho]$,
\begin{equation*}
  \exs\big(\nlev_\delta g \cap \ball_X(\rho); \nlev_{\epsilon} f\big) \leq \exs\big(\epi g \cap \ball_{X\times \reals}(\rho); \epi f\big) \leq \hatsetd_\rho(\epi f, \epi g)
\end{equation*}
provided that $\epsilon > \delta +  \exs(\epi g \cap \ball_{X\times \reals}(\rho); \epi f)$.
\end{proposition}

A parallel development is possible for set-valued mappings from a metric space $(X,d_X)$ to a metric space $(Y,d_Y)$. The values of a set-valued mapping $S:X\tto Y$ are the subsets $S(x) \subset Y$, $x\in X$, and the {\it graph} of $S$ is
\[
\gph S := \big\{(x,y)\in X \times Y~\big|~y \in S(x)\big\}.
\]
The truncated Hausdorff distance between such graphs requires a metric on $X\times Y$. Throughout, we adopt the product metric  $((x,y),(\bar x, \bar y))\mapsto \max\{d_X(x,\bar x), d_Y(y,\bar y)\}$. The centroid is likewise constructed from those of $X$ and $Y$. A prime example of such mappings is the subgradient mapping $\partial f:X\tto X$ for a convex function $f$ on a Hilbert space $X$. We recall that a function $f:X\to\Reals$ is {\it proper} if $\epi f \neq \emptyset$ and $f>-\infty$. It is {\it lower-semicontinuous} (lsc) if $\epi f$ is closed as a subset of $X\times \reals$.

\begin{proposition}{\rm (approximation of subgradient mappings \cite{AttouchWets.91}).}\label{prop:subgradmappings}
For a Hilbert space $X$, proper lsc convex functions $f,g:X\to\Reals$, and $\rho\in \preals$ exceeding $\dist(0,\epi f)$ and $\dist(0,\epi g)$, there exist $\kappa,\bar\rho\in \preals$ such that
\[
\hatsetd_\rho(\gph \partial f, \gph \partial g) \leq \kappa \sqrt{\hatsetd_{\bar \rho}(\epi f,\epi g)}.
\]
\end{proposition}
Explicit expressions for the constants $\kappa$ and $\bar\rho$ in the proposition are available in \cite{AttouchWets.91}. Section 5 establishes that $\hatsetd_\rho(\gph \partial f, \gph \partial g)$ bounds the discrepancy between near-solutions of the generalized equations $0\in \partial f(x)$ and $0\in \partial g(x)$. Thus, the proposition provides yet another way of bounding the distance between  minimizers of $f$ and those of $g$ in the convex case.

We can bring forward the effect of a  constraint set $C\subset X$ when the function of interest is expressed as $f + \iota_C$, where
\[
\iota_C(x) := 0 \mbox{ if } x\in C \mbox{ and } \iota_C(x) := \infty \mbox{ otherwise}.
\]
Then, optimality conditions can be stated using normal cones. For example, if $C\subset\reals^n$  and $f:\reals^n\to \reals$ are convex, then the generalized equation $0 \in \partial f(x) + N_C(x)$ characterizes minimizers of $f+\iota_C$, where $N_C(x)$ is the normal cone of $C$ at $x$ in the sense of convex analysis; see \cite[6.C]{VaAn}. Consequently, it becomes important to examine the graph  of a normal cone mapping $N_C:X\tto X$ and its approximations.

\begin{proposition}\label{prop:normalcone}{\rm (approximation of normal cone mappings).}
For closed convex subsets $C,D$ of a Hilbert space and $\rho\in \preals$ exceeding $\dist(0,C)$ and $\dist(0,D)$, there exist $\kappa,\bar\rho\in \preals$ such that
\[
\hatsetd_\rho(\gph N_C, \gph N_D) \leq \kappa \sqrt{\hatsetd_{\bar\rho}(C,D)}.
\]
\end{proposition}
\state Proof. In view of Cor. \ref{cor:indicator} below, the result is a direct application of Prop. \ref{prop:subgradmappings} to the functions $f = \iota_C$ and $g = \iota_D$.\eop

These preliminary facts point to a strategy for stability and error analysis of optimization and variational problems that extends much beyond the convex case: estimate the truncated Hausdorff distances between the relevant constraint sets, graphs, and/or epigraphs, which then immediately provide bounds on the discrepancy between solutions. The next sections develop practical guidelines for computing the truncated Hausdorff distance and illustrate the strategy in concrete instances.

\section{Distances between Sets}

We start with results about product sets, unions, and convex hulls. The main theorem of the section bounds the truncated Hausdorff distance between images of sets under Lipschitz continuous set-valued mappings.

\begin{proposition}{\rm (product sets).}\label{pProductSets} For each $i=1, \dots, m$, suppose that $C_i, D_i$ are subsets of a metric space $(X_i,d_{X_i})$ with centroid $x^{\rm ctr}_i$ and $X = X_1 \times \dots \times X_m$ is equipped with the metric $d_X = \max_{i=1, \dots, m} d_{X_i}$ and centroid $x^{\rm ctr}=(x^{\rm ctr}_i, \dots, x^{\rm ctr}_m)$. Then, with $C = C_1 \times \dots \times C_m$ and $D = D_1 \times \dots \times D_m$,
\[
\hatsetd_\rho(C,D) \leq \max_{i=1, \dots, m} \hatsetd_\rho(C_i, D_i) \mbox{ for any } \rho\in \preals.
\]
If $C \cap \ball_X(\rho)$ and $D\cap\ball_X(\rho)$ are nonempty, then the relation holds with equality.
\end{proposition}
\state Proof. Let $\eta = \max_{i=1, \dots, m} \hatsetd_\rho(C_i, D_i)$, $x = (x_1, \dots, x_m)\in C\cap \ball_X(\rho)$, and $\epsilon>0$. Since $x_i\in C_i\cap \ball_{X_i}(\rho)$ and $\dist(x_i,D_i)\leq \exs(C_i \cap \ball_{X_i}(\rho); D_i) \leq \eta$, there exists $y_i\in D_i$ with $d_{X_i}(x_i,y_i) \leq \eta+\epsilon$. We can repeat this construction for all $i$ and obtain $y = (y_1, \dots, y_m)$. Then, $d_X(x,y) = \max_{i=1, \dots, m} d_{X_i}(x_i,y_i) \leq \eta+\epsilon$. Thus, $\dist(x,D) \leq \eta+\epsilon$ and also $\exs(C\cap\ball_X(\rho); D) \leq \eta + \epsilon$, which holds trivially also when $C\cap\ball_X(\rho)=\emptyset$. Repeating the argument with the roles of $C$ and $D$ reversed establishes that $\hatsetd_\rho(C,D)\leq \eta + \epsilon$. Since this holds for all $\epsilon>0$, $\hatsetd_\rho(C,D)\leq \eta$ and the first conclusion holds.

To establish the inequality the other way, let $x=(x_1, \dots, x_m) \in C\cap \ball_X(\rho)$, $\epsilon>0$, and $i\in \{1, \dots, m\}$. Then, there exists $y=(y_1, \dots, y_m)\in D$ such that
\[
\dist(x_i,D_i)-\epsilon \leq d_{X_i}(x_i,y_i) - \epsilon \leq d_X(x,y) - \epsilon \leq \dist(x,D) \leq \hatsetd_\rho(C,D).
\]
Since $x\in C\cap \ball_X(\rho)$ is arbitrary, $\exs(C_i\cap \ball_{X_i}(\rho); D_i) \leq \hatsetd_\rho(C,D) + \epsilon$. A similar argument with the roles of $C$ and $D$ reversed, allows us to conclude that $\exs(D_i\cap \ball_{X_i}(\rho); C_i) \leq \hatsetd_\rho(C,D) + \epsilon$. Thus, $\hatsetd_\rho(C_i,D_i)\leq \hatsetd_\rho(C,D) + \epsilon$. Since $i$ and $\epsilon$ are arbitrary, the conclusion follows.\eop

\begin{corollary}\label{cor:indicator}{\rm (indicator functions).}
For subsets $C,D$ of a metric space and $\rho\in \preals$,
\[
\hatsetd_\rho(\epi \iota_C,\epi \iota_D) = \hatsetd_\rho(C,D).
\]
\end{corollary}
\state Proof. By Prop. \ref{pProductSets}, $\hatsetd_\rho(\epi \iota_C,\epi \iota_D)$ $=$ $\hatsetd_\rho(C\times \preals,D\times \preals)$ $=$ $\max\{\hatsetd_\rho(C,D)$, $\hatsetd_\rho(\preals,\preals)\}$ $=$ $\hatsetd_\rho(C,D)$ as long as $C \cap \ball_X(\rho)$ and $D\cap\ball_X(\rho)$ are nonempty. If one or both of these sets are empty, the corollary holds trivially.\eop

\begin{proposition}\label{prop:union}{\rm (union of sets).} For a metric space $X$, $\{C_\alpha,D_\alpha\subset X, \alpha\in A\}$, with $A$ being an arbitrary set, and $\rho\in \preals$,
\[
\hatsetd_\rho\Bigg( \bigcup_{\alpha\in A} C_\alpha, \bigcup_{\alpha\in A} D_\alpha\Bigg) \leq \sup_{\alpha\in A} \hatsetd_\rho(C_\alpha,D_\alpha).
\]
\end{proposition}
\state Proof. Let $C = \cup_{\alpha\in A} C_\alpha$, $D = \cup_{\alpha\in A} D_\alpha$, and $\eta = \sup_{\alpha\in A} \hatsetd_\rho(C_\alpha,D_\alpha)$. Suppose that $x\in C\cap \ball_X(\rho)$. Then, there exists $\alpha \in A$ such that $x\in C_{\alpha}$. Since $D_\alpha\subset D$ and $x\in C_\alpha\cap\ball_X(\rho)$,
\[
\dist(x, D)  \leq \dist(x,D_\alpha) \leq \exs(C_\alpha\cap\ball_X(\rho); D_\alpha) \leq \hatsetd_\rho(C_\alpha,D_\alpha) \leq \eta.
\]
The arbitrary choice of $x\in  C\cap \ball_X(\rho)$ allows us to conclude that $\exs(C\cap\ball_X(\rho); D) \leq \eta$. The roles of $C$ and $D$ can be reversed yielding the conclusion.\eop

There is no similar result for intersections. A revealing example is furnished already on $\reals$ by $C_1=C_2 = \{0\}$, $D_1 = \{-\epsilon\}$, and $D_2 = \{\epsilon\}$ with $\epsilon>0$. Then, $\hatsetd_\rho(C_1\cap C_2,D_1\cap D_2) = \infty$ because $D_1\cap D_2=\emptyset$. However, $\hatsetd_\rho(C_i,D_i) = 2\epsilon$ for $\rho\geq \epsilon$ and $i=1,2$.
The difficult occurs even if $C_1\cap C_2$ and $D_1\cap D_2$ have nonempty interiors. Consider $C_1 = D_1 = [-1,0]\cup [1,2]$ and $C_2 = [-1,0]\cup [2,3]$ and $D_2 = [-1,0]\cup [2+\epsilon, 3]$ with $\epsilon \in (0,1)$. Then, $C_1\cap C_2 = [-1,0]\cup \{2\}$, $D_1\cap D_2 = [-1,0]$, and $\hatsetd_\rho(C_i,D_i) \leq \epsilon$ for $i=1,2$ and $\rho \geq 3$. Still, $\hatsetd_\rho(C_1\cap C_2,D_1\cap D_2) = 2$. In the convex case, having intersections with nonempty interior remedy the situation to a large extent; see \cite[Cor. 2.5]{AzePenot.90b}. In the general case, however, it is difficulty to say more than $\exs(\cap_{\alpha\in A} C_\alpha; \cap_{\alpha\in A} D_\alpha^+) \leq 0$, where $D_\alpha^+ = \{x \in X~|~\dist(x,D_\alpha)\leq \exs(C_\alpha; D_\alpha)\}$ for $\alpha \in A$, which nevertheless provides guidance towards constructing outer approximations.

 For large enough $\rho$, the operation of taking the convex hull is non-expansive  under $\hatsetd_\rho$. We denote by $\con C$ the {\it convex hull} of a set $C$ and $\nats$ the natural numbers.

\begin{proposition}{\rm (convex hulls).}
For subsets $C$ and $D$ of a normed linear space $X$,
\[
\hatsetd_\rho(\con C, \con D) \leq \hatsetd_{\rho}(C,D)
\]
when $\rho \in [0,\infty]$ is such that $C,D\subset \ball_X(\rho)$.
\end{proposition}
\state Proof. Suppose that $x\in \con C \cap \ball_X(\rho)$. Thus, there exist $r\in\nats$, $x^1$, $\dots,$ $x^r\in C$, and $\alpha_1, \dots, \alpha_r\geq 0$, with $\sum_{i=1}^r \alpha_i = 1$ such that $x = \sum_{i=1}^r \alpha_i x^i$. Let $\epsilon>0$. Since $x^i\in C\cap\ball_X(\rho)$, there exists $y^i \in D$ with
$\|x^i - y^i\|-\epsilon \leq \dist(x^i,D) \leq \exs(C\cap\ball_X(\rho); D)\leq \hatsetd_{\rho}(C,D)$. For $y = \sum_{i=1}^r \alpha^i y^i$, $\|x-y\| \leq \sum_{i=1}^r \alpha_i\|x^i - y^i\| \leq \hatsetd_{\rho}(C,D) + \epsilon$. Thus, $\dist(x,\con D) \leq \hatsetd_{\rho}(C,D) + \epsilon$ because $y\in \con D$. Since $\epsilon$ and $x$ are arbitrary, $\exs(\con C \cap\ball_X(\rho); \con D) \leq  \hatsetd_{\rho}(C,D)$. The conclusion then follows by symmetry.\eop

The difficulty with unbounded sets and a finite $\rho$ is illustrated by $C = \{\lambda(-1,1)$, $\lambda(1,-1)\}$ $\subset\reals^2$ and $D = \{\lambda(1,1), \lambda(-1,-1)\}\subset\reals^2$, with $\lambda>0$. For the norm $\|\cdot\|_\infty$ and $\rho<\lambda$, $\hatsetd_\rho(\con C, \con D) =  \rho$ but $\hatsetd_\rho(C,D) = 0$. Near the origin $C$ and $D$ look the same (empty), but their convex hulls are locally rather different.

Next, we turn the focus towards images of sets, which provide foundations for several subsequent results.
For metric spaces $(X,d_X)$ and $(Y,d_Y)$, we say that a set-valued mapping $S:X\tto Y$ is {\it Lipschitz continuous with modulus} $\kappa:\preals\to \preals$ {\it relative to} $\rho^*\in [0,\infty]$ if
\[
\hatsetd_{\rho^*}\big(S(x),S(\bar x)\big) \leq \kappa(\rho) d_X(x,\bar x) \mbox{ for } x,\bar x\in \ball_X(\rho) \mbox{ and } \rho\in \preals.
\]
We retain this terminology also for point-valued mappings, in which case the left-hand side amounts to the truncated Hausdorff distance between two points.

The image of $C\subset X$ under a set-valued mapping $S:X\tto Y$ is the set $S(C) := \cup_{x\in C} S(x)$. The corresponding inverse set-valued mapping is $S^{-1}(y) := \{x\in X~|~y \in S(x)\}$ for $y\in Y$. Moreover, for any nonempty $C\subset X$ and $f:X\to \Reals$, $\ninf_C f := \ninf\{f(x)~|~x\in C\}$ and $\nsup_C f := \nsup\{f(x)~|~x\in C\}$. When $C$ is empty, $\ninf_C f = \infty$ and $\nsup_C f = -\infty$.

\begin{theorem}\label{thm:Lipmapping}{\rm (images under Lipschitz mappings).}
Suppose that $(X,d_X)$ and $(Y,d_Y)$ are metric spaces, with centroids $x^{\rm ctr}$ and $y^{\rm ctr}$, respectively,  $\rho\in \preals$, and $S,T:X\tto Y$ are nonempty-valued Lipschitz continuous with common modulus $\kappa:\preals \to \preals$ relative to $\rho^*\in [0,\infty]$. Then,
for any nonempty $C,D\subset X$,
\[
\hatsetd_\rho\big(S(C),T(D)\big) \leq \nsup_{x\in \ball_X(\bar \rho)}\hatsetd_{\rho^*}\big(S(x),T(x)\big) + \kappa(\hat \rho)\hatsetd_{\bar\rho}(C,D)
\]
provided that $\rho^* > 2\rho + \max\{\dist(y^{\rm ctr}, S(C)), \dist(y^{\rm ctr}, S(D)), \dist(y^{\rm ctr}, T(D))\}$, $\bar \rho > 0$ exceeds
\[
\sup_{y\in U(E) \cap \ball_Y(\rho^*)} \Big\{\inf_{U^{-1}(y) \cap E} d_X(\cdot,x^{\rm ctr})\Big\} \mbox{ for } U = S, T \mbox{ and } E = C,D,
\]
and $\hat \rho > \bar \rho + \hatsetd_{\bar\rho}(C,D)$.
\end{theorem}
\state Proof. First, we bound $\hatsetd_{\rho^*}(S(C),S(D))$. Suppose that $\bar y \in S(C)\cap \ball_Y(\rho^*)$. Then there exists $\bar x \in S^{-1}(\bar y) \cap C$ such that $d_X(\bar x, x^{\rm ctr}) \leq \bar \rho$, i.e., $\bar x \in C \cap \ball_X(\bar \rho)$. Let $\epsilon\in (0, \hat \rho - \bar \rho - \hatsetd_{\bar \rho}(C,D))$. There exists $x\in D$ such that $\hatsetd_{\bar \rho}(C,D) \geq \exs\big(C \cap \ball_X(\bar \rho); D\big) \geq \dist(\bar x, D) \geq d_X(\bar x, x) - \epsilon$. Thus, $d_X(x, x^{\rm ctr}) \leq d_X(\bar x, x^{\rm ctr}) + d_X(\bar x, x) \leq \bar \rho +  \hatsetd_{\bar \rho}(C,D) + \epsilon \leq \hat \rho$ so that both $\bar x$ and $x$ are in $\ball_X(\hat \rho)$. There exists $y \in S(x)$ such that $d_Y(\bar y,y) \leq \dist(\bar y, S(x)) + \epsilon$, which implies that $y\in S(D)$. Then,
\begin{align*}
d_Y(\bar y, y) &\leq \dist\big(\bar y, S(x)\big) + \epsilon \leq \exs\big(S(\bar x)\cap \ball_Y(\rho^*); S(x)\big) + \epsilon\\
& \leq \hatsetd_{\rho^*}\big(S(\bar x), S(x)\big) + \epsilon \leq \kappa(\hat \rho)d_X(\bar x, x) + \epsilon \leq \kappa(\hat \rho)\hatsetd_{\bar \rho}(C,D) + (\kappa(\hat \rho)+1)\epsilon,
\end{align*}
which implies that $\exs(S(C) \cap\ball_Y(\rho^*); S(D)) \leq \kappa(\hat \rho)\hatsetd_{\bar \rho}(C,D) + (\kappa(\hat \rho)+1)\epsilon$. Repeating the arguments with the roles of $C$ and $D$ reversed and recognizing that $\epsilon$ is arbitrary, lead to
\[
\hatsetd_{\rho^*}\big(S(C),S(D)\big) \leq \kappa(\hat \rho)\hatsetd_{\bar \rho}(C,D).
\]

Second, we bound $\hatsetd_{\rho^*}\big(S(D),T(D)\big)$. Suppose that $\bar y \in S(D)\cap \ball_Y(\rho^*)$. Then there exists $\bar x \in S^{-1}(\bar y) \cap D$ such that $d_X(\bar x, x^{\rm ctr}) \leq \bar \rho$, i.e., $\bar x \in D \cap \ball_X(\bar \rho)$. Let $\epsilon>0$. There exists $y \in T(\bar x)$ such that $d_Y(\bar y,y) \leq \dist(\bar y, T(\bar x)) + \epsilon$, which implies that $y\in T(D)$. Then,
\begin{align*}
d_Y(\bar y, y) &\leq \dist\big(\bar y, T(\bar x)\big) + \epsilon \leq \exs\big(S(\bar x)\cap \ball_Y(\rho^*); T(\bar x)\big) + \epsilon\\
& \leq \hatsetd_{\rho^*}\big(S(\bar x), T(\bar x)\big) + \epsilon \leq \sup_{x\in\ball_X(\bar\rho)} \hatsetd_{\rho^*}\big(S(x),T(x)\big) + \epsilon,
\end{align*}
which implies that $\exs(S(D) \cap\ball_Y(\rho^*); T(D)) \leq \sup_{x\in\ball_X(\bar\rho)} \hatsetd_{\rho^*}\big(S(x),T(x)\big) + \epsilon$. Again by symmetry and the fact that $\epsilon$ is arbitrary, we conclude that
\[
\hatsetd_{\rho^*}\big(S(D),T(D)\big) \leq \sup_{x\in\ball_X(\bar\rho)} \hatsetd_{\rho^*}\big(S(x),T(x)\big).
\]
The result now follows by Prop. \ref{prop:triangle}.\eop

The requirement on $\bar \rho$ in the proposition is most easily verified when $C$ and $D$ are bounded, but other possibilities exist for example under a Lipschitz property on the inverse set-valued mappings. An example of this appears in Cor. \ref{cor:compAffine} below.

Sums of sets arise among other places in subdifferential calculus: For functions $f_1$ and $f_2$, the set of subgradients $\partial (f_1 + f_2)(x) = \partial f_1(x) +\partial f_2(x)$ under appropriate assumptions \cite[Sec. 10.9]{VaAn}; here and below subgradients are of the general kind\footnote{For $f:\reals^n\to \Reals$ and a point $\bar x$ where $f$ is finite, we recall that $v\in \widehat\partial f(\bar x)$ (a subgradient of the regular kind) if and only if $f(x) \geq f(\bar x) + \langle v, x-\bar x\rangle + o(\|x-\bar x\|_2)$. Moreover, $v\in \partial f(\bar x)$ (a subgradient of the general kind) if and only if
there exist $v^\nu\to v$ and $x^\nu\to x$, with $f(x^\nu)\to f(x)$, such that $v^\nu \in \widehat\partial f(x^\nu)$. In the convex case, regular and general subgradients coincide.} \cite{VaAn,Mordukhovich.13b}. Of course, the previous theorem could be used to establish a result about sums. We pursue a direct approach, with a proof in the appendix, as it is instructive and also brings forth a possible adjustment in the case of unbounded sets.

\begin{proposition}\label{prop:sumofsets}{\rm (sums of sets).}
For a normed linear space $X$, nonempty sets $\{C_i,D_i\subset X, i=1, \dots, m\}$, and $\rho\in \preals$,
\[
\hatsetd_\rho\Bigg( \sum_{i=1}^m C_i, \sum_{i=1}^m D_i\Bigg) \leq \sum_{i=1}^m \hatsetd_{\rho}(C_i,D_i)
\]
provided that $C_i,D_i \subset \ball_X(\rho)$ for all $i=1, 2, \dots, m$. If $C_i,D_i \subset \ball_X(\rho)$ holds only for $i = 2, 3, \dots, m$, then the inequality remains valid as long as $\hatsetd_{\rho}(C_1,D_1)$ is replaced by $\hatsetd_{m\rho}(C_1,D_1)$.
\end{proposition}

A motivation for allowing one unbounded set merges when studying a locally Lipschitz continuous function $f:\reals^n\to \reals$, a nonempty closed set $C\subset\reals^n$, and the optimality condition  $0\in \partial f(x)+N_C(x)$ \cite[Exer 10.10]{VaAn}, where $N_C(x)$ is the normal cone of $C$ at $x$ in the general sense \cite{VaAn,Mordukhovich.13b}, i.e., $N_C(x) = \partial \iota_C(x)$. Here, $\partial f(x)$ is bounded, but $N_C(x)$ is not in the interesting cases. We observe that if there are two or more unbounded sets, then the assertion in the proposition fails. For an example in $\reals^2$, let $C_1 = \{\lambda (1, 1+\delta)~|~\lambda\geq 0\}$, $C_2 = \{\lambda(-1,-1+\delta)~|~\lambda\geq 0\}$, with $\delta>0$, $D_1 =\{\lambda(1,1)~|~\lambda\geq 0\}$, and $D_2 = \{\lambda(-1,-1)~|~\lambda \geq 0\}$. All the sets are rays and therefore unbounded. Now, $\hatsetd_\rho(C_i,D_i) \leq \delta\rho$ for $i=1, 2$. However, because $C_1+C_2$ is ``nearly'' the halfspace $\{(x_1,x_2)~|~x_1-x_2\leq 0\}$ for small $\delta$ but $D_1 + D_2 = \{(x_1,x_2)~|~x_1 = x_2\}$, $\hatsetd_\rho(C_1+C_2, D_1+D_2)=\rho$.

The inequality in the proposition is sharp because for $x,y,z\in X$ and $C_1 = \{x\}$, $C_2 = \{y\}$, $D_1 = \{x+z\}$, and $D_2 = \{y+z\}$, we have $\hatsetd_\rho(C_1+C_2,D_1+D_2) = 2\|z\|$ and $\hatsetd_\rho(C_i,D_i)=\|z\|$ for $i=1,2$ for sufficiently large $\rho$. Still, we can have strict inequality. For example, $x,y\in X$, $x\neq y\neq 0$, and $C_1 = \{x\}$, $C_2 = \{-x\}$, $D_1 = \{y\}$, and $D_2 = \{-y\}$, we have $\hatsetd_\rho(C_1+C_2,D_1+D_2) = 0$ and $\hatsetd_\rho(C_i,D_i)=\|x-y\|$ for $i=1,2$ for sufficiently large $\rho$.

\begin{corollary}\label{prop:setmultiplication}{\rm (set multiplications).}
For nonempty subsets $C$ and $D$ of a normed linear space, nonzero $\lambda,\mu\in\reals$, and $\rho\in\preals$,
\[
\hatsetd_\rho(\lambda C,\mu D) \leq \bar\rho|\lambda-\mu| + \max\{|\lambda|,|\mu|\}\hatsetd_{\bar\rho}(C,D),
\]
when $\bar\rho > (2\rho + \max\{|\lambda|\dist(0,C), |\lambda|\dist(0,D), |\mu|\dist(0,D)\})\max\{|\lambda^{-1}|,|\mu^{-1}|\}$.
\end{corollary}
\state Proof. The result follows from Thm. \ref{thm:Lipmapping} by setting $S(x) = \lambda x$ and $T(x) = \mu x$.\eop

We end the section by recording a useful fact about the distance between level-sets of two convex functions, which extends \cite[Prop. 7.68]{VaAn} by allowing the functions to be different.

\begin{proposition}\label{prop:convexlevel}{\rm (level-sets; convex case).}
For $\rho \in \preals$, $\alpha,\beta \in [-\rho,\rho]$, and proper convex lsc functions $f,g:\reals^n\to \Reals$, suppose that $\alpha>\inf f$, $\beta>\inf g$, $\nargmin f\neq\emptyset$, and $\nargmin g\neq \emptyset$. Then, with $\eta = \hatsetd_\rho(\epi f, \epi g)$,
\[
\hatsetd_\rho(\nlev_\alpha f, \nlev_\beta g) \leq \eta + (\rho^*+\rho_0)\max\left\{\frac{\alpha + \eta - \beta}{\alpha + \eta - \inf g}, \frac{\beta + \eta - \alpha}{\beta + \eta - \inf f}\right\}
\]
provided that $\rho_0 \geq \max\{\dist(0, \nargmin f), \dist(0,\nargmin g)\}$ and $\rho^* \geq \max\{\rho_0, \rho + \hatsetd_\rho(\epi f, \epi g)\}$.
\end{proposition}
\state Proof. By Prop. 4.5 in \cite{Royset.18}, $\exs(\nlev_\alpha f \cap \ball_{\reals^n}(\rho); \nlev_{\alpha+\eta} g) \leq \eta$.  An application of Prop. 7.68 in \cite{VaAn} yields
\[
\exs\big(\nlev_{\alpha+\eta} g\cap \ball_{\reals^n}(\rho^*); \nlev_\beta g\big) \leq \frac{\alpha + \eta - \beta}{\alpha + \eta - \inf g}(\rho^* + \rho_0)
\]
whenever $\alpha + \eta > \beta$. If $\alpha + \eta \leq \beta$, then $\exs(\nlev_{\alpha+\eta} g\cap \ball_{\reals^n}(\rho^*); \nlev_\beta g) = 0$.
Let $x\in \nlev_\alpha f \cap \ball_{\reals^n}(\rho)$. There exists $y\in \nlev_{\alpha+\eta} g$ with $\|y - x\|\leq \eta$ so that $y\in \ball_{\reals^n}(\rho^*)$. Thus, we have established that
\[
\exs\big(\nlev_{\alpha} f\cap \ball_{\reals^n}(\rho); \nlev_\beta g\big) \leq \eta + \max\left\{0, \frac{\alpha + \eta - \beta}{\alpha + \eta - \inf g}(\rho^* + \rho_0)\right\}.
\]
Repeating the argument with the roles of $f$ and $g$ reversed leads to the conclusion.\eop

The proposition relies heavily on the assumption that $\nlev_\alpha f$ and $\nlev_\beta g$ have nonempty interiors. The next section dispenses of that requirement as well as convexity.

\section{Distances between Epigraphs of Functions}

As special sets, epigraphs offer several possibilities to specialize the results of the previous section and also develop new ones. First, we examine the Kenmochi conditions and their numerous applications including in the analysis of constrained problems with feasible sets that lack interiors. Second, we develop a series of calculus rules relying, in part, on Section 3.

For a metric space $(X,d_X)$, let the closed balls at $x\in X$ be denoted by
\[
\ball_X(x, \rho) := \{\bar x\in X~|~d_X(x,\bar x) \leq \rho\} \mbox{ for } \rho\geq 0.
\]

\subsection{Kenmochi Conditions and Applications}

An alternative expression for the truncated Hausdorff distance between epigraphs is provided by the {\it Kenmochi conditions}, which can be traced back to \cite{Kenmochi.74}; see also \cite{AttouchWets.91}.
The following result generalizes \cite[Prop. 3.2]{Royset.18} by relaxing a lsc assumption and  establishing that the conditions provide tight estimates. A proof is provided in the appendix.

\begin{proposition}\label{prop:altform}{\rm (Kenmochi conditions). }
For a metric space $X$, functions $f,g:X\to \Reals$, both with nonempty epigraphs, and $\rho\in \preals$,
\begin{align*}
\hatsetd_\rho(\epi f,\epi g) = \ninf\Big\{ \eta\geq 0 \Big| &\ninf_{\ball_X(x,\eta)} g \leq \max\{f(x), -\rho\} + \eta, \forall x\in \nlev_\rho f \cap \ball_X(\rho)\\
                           & \ninf_{\ball_X(x,\eta)} f \leq \max\{g(x), -\rho\} + \eta, \forall x\in \nlev_\rho g \cap \ball_X(\rho) \Big\}.
\end{align*}
\end{proposition}

For $\alpha\in (0,\infty)$, a function $f:X\to \Reals$ is $\alpha$-{\it H\"{o}lder continuous with modulus} $\kappa:\preals\to \preals$ if
\[
|f(x)-f(\bar x)| \leq \kappa(\rho) \big[d_X(x,\bar x)\big]^\alpha \mbox{ for } x,\bar x\in \ball_X(\rho) \mbox{ and } \rho\in \preals.
\]
The function is {\it Lipschitz continuous with modulus} $\kappa:\preals\to\preals$ if the relation holds with $\alpha = 1$.

The truncated Hausdorff distance between epigraphs of functions of this kind can be bounded by an expression involving the worst pointwise difference between the functions over a set.

\begin{proposition}\label{prop:fromsupnorm}{\rm (estimates from sup-norm).} For a metric space $X$, functions $f,g:X\to \Reals$ with nonempty epigraphs, and $\rho\in \preals$, we have that
\[
\hatsetd_\rho(\epi f, \epi g) \leq \nsup_{A_\rho}|f-g|,
\]
where $A_\rho = \nlev_\rho f \cup \nlev_\rho g \cap \ball_X(\rho)$. (Supremum over an empty set is interpreted as zero in this case.) 
Suppose also that  $f$ and $g$ are $\alpha$-H\"{o}lder continuous with common modulus $\kappa:\preals\to \preals$ and $\alpha\in (0,\infty)$. Then, for any nonempty $C\subset X$,
\[
\hatsetd_\rho(\epi f, \epi g) \leq \max\big\{\exs(A_\rho; C), ~\kappa(\hat\rho)[\exs(A_\rho; C)]^\alpha + \nsup_C|f-g|\big\},
\]
provided that $\hat\rho > \rho+\exs(A_\rho; C)$.
\end{proposition}
\state Proof. The first assertion holds via Prop. \ref{prop:altform}. For the second assertion, set $\eta = \exs(A_\rho; C)$ and let $\epsilon \in (0,\hat\rho-\rho-\eta)$. Suppose that $x\in \nlev_\rho f \cap \ball_X(\rho)$. Then, there exists $\bar x\in C$ with $d_X(x,\bar x) \leq \bar\eta = \eta+\epsilon$ and
\[
\ninf_{\ball_X(x,\bar\eta)} g \leq g(\bar x) \leq f(\bar x) + \nsup_C|f-g|\leq \max\{f(x), -\rho\} + \kappa(\hat\rho)\bar\eta^\alpha + \nsup_C|f-g|.
\]
A similar result holds with the roles of $f$ and $g$ reversed. Thus, by Prop., \ref{prop:altform} $\hatsetd_\rho(\epi f, \epi g)$ $\leq$ $\max\{\bar\eta, \kappa(\hat\rho)\bar\eta^\alpha + \nsup_C|f-g|\}$. Since $\epsilon$ is arbitrary, $\bar \eta$ can be replaced by $\eta$ and the second conclusion holds.\eop

{\it Example 1: sample average approximations}. In stochastic optimization and statistical learning, $f:X\to \reals$ is often given as $f(x) = \mathbb{E}[\psi(\boldmath{\xi},x)]$, where $\psi:\Xi\times X\to \reals$ and $\mathbb{E}$ denotes the expectation under the distribution of the random vector $\boldmath{\xi}$ with values in $\Xi$. Under standard assumptions (see \cite[Ch. 14]{VaAn}, \cite[Ch. 7]{ShapiroDentchevaRuszczynski.14}), $f$ is well defined and Lipschitz continuous with modulus $\kappa:\reals_+\to \reals_+$. An approximation of $f$ could be the sample average function $f^\nu:X\to \reals$ given by $f^\nu(x) = \nu^{-1} \sum_{i=1}^\nu \psi(\xi^i, x)$, where $\xi^1, \dots, \xi^\nu\in \Xi$ are given data. Under related assumptions, $f^\nu$ is also Lipschitz continuous with the same modulus as $f$.
When $X$ is finitely compact\footnote{Recall that a metric space is finitely compact if all its balls are compact.}, $A_\rho$ in Prop. \ref{prop:fromsupnorm} is compact and it is possible to construct for any $\epsilon>0$ a set $C$ consisting of only a finite number of points and still have $\exs(A_\rho; C)\leq\epsilon$. Since $C$ is finite, there exists a variety of ways of bounding $\nsup_C|f-f^\nu|$, say by $\delta$, using the theory of large deviations; see for example \cite[Ch. 7]{ShapiroDentchevaRuszczynski.14}. Prop. \ref{prop:fromsupnorm} then gives that $\hatsetd_\rho(\epi f, \epi f^\nu) \leq \max\{\epsilon, \kappa(\hat\rho)\epsilon + \delta\}$ when $\hat\rho > \rho+\epsilon$.\\

The next result extends \cite[Prop. 3.3]{Royset.18} by moving from indicator functions to general functions and from Lipschitz to H\"{o}lder continuous functions; see also \cite{AttouchWets.91,AzePenot.90b} for results on sums in the convex case.

\begin{proposition}\label{prop:lip}{\rm (sums under H\"{o}lder continuity).}
For a metric space $X$, functions $f_i,g_i:X\to\Reals$, $i=1, 2$, where $f_1, g_1$ are $\alpha$-H\"{o}lder continuous with common modulus $\kappa:\preals\to \preals$, $\alpha\in (0,\infty)$, and both $\epi(f_1+f_2)$ and $\epi(g_1+g_2)$ are nonempty. Then, for $\rho\in \preals$, 
\[
\hatsetd_\rho\big(\epi (f_1+f_2), \epi (g_1+g_2) \big) \leq \nsup_{A_\rho} |f_1-g_1| + \eta + \kappa(\hat\rho)\eta^\alpha
\]
where $\eta = \hatsetd_{\bar\rho}(\epi f_2,\epi g_2)$, provided that $A_\rho = \nlev_\rho (f_1+f_2) \cup \nlev_\rho (g_1 + g_2) \cap \ball_X(\rho)\neq\emptyset$, $\bar\rho \geq \rho + \max\{\nsup_{\ball_X(\rho)} |f_1|, \nsup_{\ball_X(\rho)} |g_1|\}$, and $\hat\rho > \rho + \eta$.
\end{proposition}
\state Proof. Let $\epsilon \in (0, \hat\rho -\rho - \eta)$ and $x\in \nlev_\rho (f_1+ f_2) \cap \ball_X(\rho)$. Then, $f_2(x) \leq \rho -f_1(x) \leq \bar\rho$. First, suppose that $f_2(x) \geq -\bar \rho$ so that $(x,f_2(x)) \in \epi f_2 \cap \ball_{X\times\reals}(\bar \rho)$. Consequently, there is $(\bar x, \bar \alpha)$ $\in$ $\epi g_2$ with $d_X(x,\bar x) \leq \eta+\epsilon$ and $|\bar\alpha - f_2(x)|\leq \eta+ \epsilon$. Thus, $g_2(\bar x) \leq \bar \alpha\leq f_2(x) + \eta + \epsilon$ and
\begin{align*}
  \ninf_{\ball_X(x,\eta+\epsilon)} g_1 + g_2 &\leq g_1(\bar x) + g_2(\bar x) = g_1(\bar x) - g_1(x) + g_1(x) - f_1(x) + f_1(x) + g_2(\bar x)\\
  & \leq \kappa(\hat\rho)(\eta+\epsilon)^\alpha + \nsup_{A_\rho}|f_1 - g_1| + f_1(x) + f_2(x) + \eta + \epsilon\\
  & \leq \max\big\{f_1(x) + f_2(x), -\rho\big\} + \nsup_{A_\rho}|f_1 - g_1| +  \kappa(\hat\rho)(\eta+\epsilon)^\alpha + \eta + \epsilon.
\end{align*}
Second, suppose that $f_2(x) < -\bar \rho$. Then, $(x,-\bar \rho) \in \epi f_2 \cap \ball_{X\times\reals}(\bar \rho)$ and there is $(\bar x, \bar \alpha) \in \epi g_2$ with $d_X(x,\bar x) \leq \eta+\epsilon$ and $|\bar\alpha +\bar \rho|\leq \eta+ \epsilon$. Thus, $g_2(\bar x) \leq \bar \alpha \leq -\bar\rho + \eta+\epsilon$ and, similar to above,
\begin{align*}
  \ninf_{\ball_X(x,\eta+\epsilon)} g_1 + g_2 &\leq  g_1(\bar x) - g_1(x) + g_1(x) - f_1(x) + f_1(x) + g_2(\bar x)\\
  & \leq \kappa(\hat\rho)(\eta+\epsilon)^\alpha + \nsup_{A_\rho}|f_1 - g_1| + f_1(x) -\bar\rho + \eta + \epsilon\\
  & \leq \max\big\{f_1(x) + f_2(x), -\rho\big\} + \nsup_{A_\rho}|f_1 - g_1| +  \kappa(\hat\rho)(\eta+\epsilon)^\alpha + \eta + \epsilon.
\end{align*}
The last inequality follows because $f_1(x) - \bar\rho \leq \nsup_{\ball_X(\rho)}|f_1| - \bar\rho \leq -\rho$. Thus, in both cases, we obtain the same upper bound on $\ninf_{\ball_X(x,\eta+\epsilon)} g_1 + g_2$. Repeating these arguments with the roles of $f_1,f_2$ switched with those of $g_1,g_2$, we obtain via Prop. \ref{prop:altform} that $\hatsetd_\rho\big(\epi (f_1+f_2), \epi (g_1+g_2)\big) \leq \max\{\eta + \epsilon, \nsup_{A_\rho}|f_1 - g_1| +  \kappa(\hat\rho)(\eta+\epsilon)^\alpha + \eta+ \epsilon\}$. Since $\epsilon$ is arbitrary, the conclusion follows.\eop

{\it Example 1: continued}. Suppose that in addition to $f$ the problem of interest involves a ``regularizer'' $r:X\to [0,\infty)$, which is common in statistical learning, i.e., we aim to minimize $f + r$. We may want to examine the stability of solutions under changes to $r$. Let $r^\nu:X\to [0,\infty)$ be such an alternative regularizer. A prime example is when $r = 0$ and we want to quantify the effect of the regularizer $r^\nu$. We are therefore interested in comparing $\epi (f + r)$ to $\epi (f^\nu + r^\nu)$. Suppose that $r$ and $r^\nu$ are $\alpha$-H\"{o}lder continuous with common modulus $\mu:\preals\to \preals$ and $\alpha\in (0,\infty)$, and $X = \reals^n$. A possible choice is to have $r^\nu(x) = \sum_{j=1}^n s^\nu(x_j)$ with $s^\nu(\tau) = \lambda |\tau| - \nu\tau^2/2$ when $|\tau|\leq \lambda/\nu$ and $s^\nu(\tau) = \lambda^2/(2\nu)$ otherwise, with $\lambda>0$ being a parameter. This makes $r^\nu$ a nonconvex function with Lipschitz modulus $\lambda$ globally. An even more aggressive regularizer would be $s^\nu(\tau) = \nu^{-1}\sqrt{|\tau|}$, possibly further scaled, which is nonconvex but $1/2$-H\"{o}lder continuous. Regardless, Prop. \ref{prop:lip} establishes that
\[
\hatsetd_\rho\big(\epi (f+r), \epi (f^\nu+r^\nu) \big) \leq \nsup_{A_\rho} |r-r^\nu| + \eta + \mu(\hat\rho)\eta^\alpha
\]
where $\eta = \hatsetd_{\bar\rho}(\epi f,\epi f^\nu)$ can be expressed in terms of $\kappa$, $\epsilon$, and $\delta$, and $A_\rho$ and $\hat\rho$ are sufficiently large as stipulated by the proposition. In particular when $r=0$, this error bound provides guidance on how fast the regularizer should vanish as the sample size $\nu$ grows. Typically, the sample error $\delta$ is of order $\nu^{-1/2}$, which indicates that $r^\nu$ should vanish at the same rate at least when $\alpha = 1$.\\

{\it Example 2: disjunctive programming}.  Suppose that $\{C_\alpha, \alpha\in A\}$ is a collection of nonempty subsets of a Hilbert space $X$ and $c\in X$. Disjunctive programming studies problems of the form minimize $\langle c, x\rangle$ subject to $x\in \cup_{\alpha \in A} C_\alpha$. The effect of replacing $c$ by $d\in X$ and the sets by $\{D_\alpha \neq \emptyset, \alpha\in A\}$ on the minimum value and set of near-minimizers can be bounded by Prop. \ref{prop:solestimates} via Prop. \ref{prop:lip} and Prop. \ref{prop:union}. Specifically, let $f(x) = \langle c, x\rangle$ if $x\in C = \cup_{\alpha \in A} C_\alpha$ and $f(x) = \infty$ otherwise. Likewise, $g(x) = \langle d, x\rangle$ if $x\in D = \cup_{\alpha \in A} D_\alpha$ and $g(x) = \infty$ otherwise. Since $\inf_{x\in\ball_X(\rho)} \langle c, x\rangle \geq -\rho\|c\|$ and similarly with $c$ replaced by $d$, $\bar\rho$ can be set to $\rho(1 + \max\{\|c\|,\|d\|\})$ in Prop. \ref{prop:lip} and, in view of the Lipschitz continuity of $\langle c, \cdot\rangle$ and $\langle d, \cdot\rangle$,
\begin{align*}
\hatsetd_\rho(\epi f, \epi g) &\leq \rho\|c-d\| + \big(1+\max\{\|c\|,\|d\|\}\big)\hatsetd_{\bar\rho}(\epi \iota_C, \epi \iota_D)\\
& \leq \rho\|c-d\| + \big(1+\max\{\|c\|,\|d\|\}\big)\nsup_{\alpha\in A} \hatsetd_{\bar\rho}(C_\alpha, D_\alpha),
\end{align*}
where the last inequality follows by  Cor. \ref{cor:indicator} and Prop. \ref{prop:union}. Consequently, solutions of disjunctive programs exhibit a Lipschitz property in this sense under a remarkable absence of assumptions.\\

As already discussed in Section 3, intersections of sets are generally not stable under perturbations of the individual sets.  This fact is the source of many difficulties in constrained optimization. In particular, if the problem of minimizing $f_0(x)$ subject to $x\in C_\alpha$ for all $\alpha\in A$ is ``approximated'' by minimizing $g_0(x)$ subject to $x\in D_\alpha$ for all $\alpha\in A$, with both $\nsup_X|f_0 - g_0|$ and $\hatsetd_\rho(C_\alpha, D_\alpha)$ being ``small'' for all $\alpha\in A$, then their solutions can still be arbitrarily far apart.  The issue surfaces even in one dimension: for example, set $f_0(x) = g_0(x) = x$, $C_1 = D_1 = \{0, 1\}$, $C_2 = [0, 1-\epsilon]$,  and $D_2 = [\epsilon, 1]$ for $\epsilon \in (0,1)$.  Thus, a major challenge is to construct approximating problems that are associated with small truncated Hausdorff distances to their original counterparts.
We observe that in the convex case having an intersection of constraint sets with nonempty interior suffices to avoid this difficulty as long as the approximations are sufficiently accurate; see \cite[Cor. 2.5]{AzePenot.90b}.

We illustrate three cases, while neither making assumptions about the feasible sets having an interior nor being convex. Moreover, the approximations can be arbitrarily poor, i.e., we are not only considering small perturbations. This forces us to construct approximating problems that are rather {\it different} than the actual problems because simply replacing objective functions and constraint sets by approximating counterparts usually fail to achieve small solution errors as the trivial example in the previous paragraph highlights.\\

\noindent {\bf Case I.} The first case analyzes the feasibility problem of finding an $x\in \cap_{i=1}^m C_i$ when we only have approximating sets $D_1, \dots, D_m$. We construct an approximating optimization problem in a higher-dimensional space that furnishes an approximating solution of the actual feasibility problem and is computationally attractive as it ``nearly'' decomposes into $m$ subproblems.

\begin{theorem}\label{thm:feasproblem}{\rm (approximation of feasibility problem).} For subsets $C_1, \dots, C_m$ and $D_1, \dots, D_m$ of a metric space $(X,d_X)$, with centroid $x^{\rm ctr}$, $\lambda\in (0,\infty)$, $\rho> 2\lambda(m-1)\max_{i=1, \dots, m} d_X(x^{\rm ctr}, D_i)$, with $\cap_{i=1}^m C_i \cap \ball_X(\rho) \neq\emptyset$, and $\bar \rho \in (3\rho, \infty)$, suppose that the following constraint qualification holds: there exists a nondecreasing function $\psi:\preals\to \preals$ such that
\[
\dist(x_1, \cap_{i=1}^m C_i) \leq \psi\Big(\sum_{i=1}^m d_X(x_i,x_1)\Big) \mbox{ for all } x_i \in C_i \cap\ball_X(\bar \rho), ~i=1, \dots, m.
\]
Then, any solution
\[
(\bar x_1, \dots, \bar x_m) \in \nargmin\Bigg\{\lambda \sum_{i=1}^m d_X(x_i, x_1) ~\Big|~ x_i \in D_i, ~i=1, \dots, m\Bigg\} \bigcap \ball_{X^m}(\rho)
\]
satisfies
\[
\dist\Bigg(\bar x_1, ~\bigcap_{i=1}^m C_i\Bigg) \leq \frac{\bar\rho}{\lambda} + \psi\Big(\frac{\bar\rho}{\lambda}\Big) + (1+2m\lambda)\max_{i=1, \dots, m} \hatsetd_{\bar\rho}(C_i,D_i).
\]
\end{theorem}
\state Proof. Let $C = C_1\times \dots \times C_m\subset X^m$, $D = D_1\times \dots \times D_m\subset X^m$, and define $f,f^\lambda,g^\lambda:X^m\to \Reals$ to have $f(x_1, \dots, x_m) = 0$ if $(x_1, \dots, x_m)\in C$ and $x_i = x_1$ for all $i$, $f^\lambda(x_1, \dots, x_m) = \lambda \sum_{i=1}^m d_X(x_i,x_1)$ if $(x_1, \dots, x_m)\in C$, and $g^\lambda(x_1, \dots, x_m) = \lambda \sum_{i=1}^m d_X(x_i,x_1)$ if $(x_1, \dots, x_m)\in D$. Otherwise, the functions take the value $\infty$.

First, we examine the Kenmochi conditions for  $f$ and $f^\lambda$. Suppose $(x_1, \dots, x_m)$ $\in$ $\nlev_{\bar\rho} f \cap \ball_{X^m}(\bar\rho)$. (Note that $X^m = X\times \dots \times X$ is equipped with the product metric.) Then, $(x_1, \dots, x_m) \in C$ and $x_i = x_1$ for all $i$. Thus, $\ninf_{\ball_{X^m}((x_1, \dots, x_m),0)} f^\lambda \leq f^\lambda(x_1, \dots, x_m) = 0 = f(x_1, \dots, x_m)$ and the first set of Kenmochi conditions holds with $\eta = 0$. Next, suppose that $(x_1, \dots, x_m) \in \nlev_{\bar\rho} f^\lambda \cap \ball_{X^m}(\bar\rho)$. Then, $x_i\in C_i$ for all $i$ and $\lambda \sum_{i=1}^m d_X(x_i,x_1)\leq \bar\rho$. In view of the constraint qualification, this implies that
\[
\dist(x_1, \cap_{i=1}^m C_i) \leq \psi\Big(\sum_{i=1}^m d_X(x_i,x_1)\Big) \leq \psi(\bar\rho/\lambda).
\]
Let $\epsilon>0$. There exists $\bar x\in \cap_{i=1}^m C_i$ such that $\dist(x_1,\cap_{i=1}^m C_i) \geq d_X(x_1,\bar x) - \epsilon$. Certainly,
\[
d_X(x_i,\bar x) \leq d_X(x_i,x_1) + d_X(x_1,\bar x) \leq \bar\rho/\lambda + \psi(\bar\rho/\lambda)+\epsilon.
\]
Then, with $\eta = \bar\rho/\lambda + \psi(\bar\rho/\lambda)+\epsilon$,
\[
\ninf_{\ball_{X^m}((x_1, \dots, x_m),\eta)} f \leq f(\bar x, \dots, \bar x) =  0 \leq  f^\lambda(x_1, \dots, x_m)
\]
and the second set of Kenmochi conditions holds with this $\eta$.
Since $\epsilon$ is arbitrary, we have established via Prop. \ref{prop:altform} that
\[
\hatsetd_{\bar\rho}(\epi f, \epi f^\lambda) \leq \bar\rho/\lambda + \psi(\bar\rho/\lambda).
\]

Second, we estimate $\hatsetd_{\bar\rho}(\epi f^\lambda, \epi g^\lambda)$. The Lipschitz modulus of the function $(x_1, \dots, x_m)\mapsto \lambda \sum_{i=1}^m d_X(x_i,x_1)$ is the constant $2m\lambda$. By Prop. \ref{pProductSets}, Prop. \ref{prop:lip}, and Cor. \ref{cor:indicator},
\[
\hatsetd_{\bar\rho}(\epi f^\lambda, \epi g^\lambda) \leq (1+2m\lambda)\hatsetd_{\bar \rho}(C,D) \leq (1+2m\lambda)\nmax_{i=1, \dots, m} \hatsetd_{\bar\rho}(C_i,D_i).
\]

For any $\epsilon>0$, we have that
\begin{align*}
\dist\big((x^{\rm ctr}, \dots, x^{\rm ctr}, 0), \epi f\big) &\leq  \dist(x^{\rm ctr}, \cap_{i=1}^m C_i ) \leq \rho\\
\dist\big((x^{\rm ctr}, \dots, x^{\rm ctr}, 0), \epi f^\lambda\big) &\leq  \dist(x^{\rm ctr}, \cap_{i=1}^m C_i ) \leq \rho\\
\dist\big((x^{\rm ctr}, \dots, x^{\rm ctr}, 0), \epi g^\lambda\big) &\leq  2\lambda(m-1)\max_{i=1, \dots, m} d_X(x^{\rm ctr}, D_i) + \epsilon < \rho + \epsilon.
\end{align*}
Thus, $\bar\rho > 3\rho$ is sufficiently large for use in Prop. \ref{prop:triangle} and
\[
\hatsetd_\rho(\epi f, \epi g^\lambda) \leq \eta = \bar\rho/\lambda + \psi(\bar\rho/\lambda) + (1+2m\lambda)\nmax_{i=1, \dots, m} \hatsetd_{\bar \rho}(C_i,D_i).
\]
We next apply Prop. \ref{prop:solestimates} to the functions $f$ and $g^\lambda$. The conditions of the proposition is easily verified. In particular, for $(x_1, \dots, x_m)\in D$,
\[
\inf g^\lambda \leq \lambda \sum_{i=1}^m d_X(x_i,x_1) \leq \lambda (m-1) \max_{i=1, \dots, m} d_X(x_i,x_1),
\]
which together with the fact that $d_X(x_i,x_1) \leq 2\max_{i=1, \dots, m} \dist(x^{\rm ctr}, D_i) + \epsilon$ for any $\epsilon>0$ ensure that
\[
\inf g^\lambda \leq 2\lambda(m-1)\max_{i=1, \dots, m} d_X(x^{\rm ctr}, D_i) + \lambda(m-1)\epsilon.
\]
Consequently, Prop. \ref{prop:solestimates} yields $\exs\big(\nargmin g^\lambda \cap \ball_{X^m}(\rho); ~\delta\mbox{-}\nargmin f\big) \leq \eta$ for $\delta> 2\eta$. Since $\delta\mbox{-}\nargmin f = \{(x_1, \dots, x_m)\in C~|~x_i = x_1, i=1, \dots, m\}$ for $\delta\geq 0$, the conclusion holds.\eop

The constraint qualification quantifies how close the points $\{x_i\in C_i, i=1, \dots, m\}$ will be to $\cap_{i=1}^m C_i$ when the points are close to each other. An example similar to the one discussed prior to the theorem is furnished by $C_1 = D_1 = \{0,1\}$, $C_2 = [0, 1-\delta]$, with $\delta\in (0,1)$, and $D_2 = [\epsilon, 1-\delta]$, with $\epsilon \in (0,1-\delta]$, where $\hatsetd_\rho(C_i,D_i) \leq \epsilon$ for $i=1,2$ and $\rho\geq \epsilon$. Thus, $C_1 \cap C_2 = \{0\}$, but $D_1\cap D_2 = \emptyset$ and it would be futile to attempt to find a feasible point in $C_1\cap C_2$ by solving $x\in D_1\cap D_2$. However, the approximating problem of the theorem produces the desired result. Specifically, in this case we can take $\psi(\gamma) = \gamma/\delta$ for $\gamma \geq 0$. Thus, the approximating problem produces a solution with error of at most $\bar\rho(\lambda^{-1} + \delta^{-1} \lambda^{-1}) + (1+4\lambda)\epsilon$.  As $\epsilon\downto 0$, this error  vanishes as long as $\lambda$ is set appropriately, for example to $\epsilon^{-1/2}$.

In general, the rate of convergence depends on the conditioning function $\psi$. Poor conditioning requires a large $\lambda$ that in turn increases the third term in the conclusion of Thm. \ref{thm:feasproblem}.  Even in the convex case, the conditioning can be arbitrarily poor: let $C_1 = \{x\in\reals^2~|~x_2\leq 0\}$ and $C_2 = \{x\in\reals^2~|~x_1^\alpha \leq x_2\}$ for $\alpha>1$, with $C_1\cap C_2 = \{0\}$. Then, $\psi(\gamma) = \gamma^{1/\alpha}$ and $x_1\in C_1$ and $x_2\in C_2$ can be close even though $x_1$ is far from the origin for large $\alpha$. Further details about constraint qualifications arise in the following two theorems for the case of inequality constraints.\\

\noindent {\bf Case II.} The second case considers the optimization problem
\begin{equation}\label{eqn:actualproblem}
\nnmin_{x\in X} f_0(x) \mbox{ subject to } f_i(x) \leq 0 \mbox{ for } i=1, \dots, m
\end{equation}
for which the actual functions need to be approximated by $g_0, \dots, g_m$. As already mentioned, an ``approximating'' problem obtained by simply replacing $f_i$ by $g_i$ for $i= 0, 1, \dots, m$ might fail to be epigraphically close to the actual problem \eqref{eqn:actualproblem} even though $\max_{i=0, \dots, m} \nsup_{x\in X} |f_i(x) - g_i(x)|$ is small. In particular, $\{x\in X~|~g_i(x)\leq 0, i=1, \dots, m\}$ could be empty while the actual feasible set is nonempty. As an alternative, we examine for $\lambda>0$ the approximating problem
\[
\nnmin_{x\in X, y\in \reals^m} g_0(x) + \lambda \sum_{i=1}^m y_i \mbox{ subject to } g_i(x) \leq y_i, ~~y_i\geq 0 \mbox{ for } i=1, \dots, m,
\]
with variable $y = (y_1, \dots, y_m) \in \reals^m$. We see next that this approximating problem furnishes approximating solutions for \eqref{eqn:actualproblem} via Prop. \ref{prop:solestimates}.

\begin{theorem}\label{thm:slack}{\rm (approximation by constraint softening).} For a metric space $X$ and $f_i, g_i:X\to \reals$, $i=0, 1, \dots, m$, where $f_0$ and $g_0$ are Lipschitz continuous with common modulus $\kappa:\preals\to \preals$, consider the functions $f,g^\lambda:X\times \reals^m\to \Reals$ defined by
\[
f(x,y) = \begin{cases}
f_0(x) &\mbox{ if } f_i(x) \leq 0 \mbox{ and } y_i = 0 \mbox{ for all } i = 1, \dots, m\\
\infty &\mbox{ otherwise}
\end{cases}
\]
and, with $\lambda\in (0,\infty)$,
\[
g^\lambda(x,y) = \begin{cases}
g_0(x) + \lambda \sum_{i=1}^m y_i &\mbox{ if } g_i(x) \leq y_i, ~y_i \geq 0\mbox{ for all } i=1, \dots, m\\
\infty &\mbox{ otherwise.}
\end{cases}
\]
Then\footnote{Here we use the product metric on $X\times \reals^m$ constructed from the sup-norm on $\reals^m$.}, for $\rho\in \preals$,
\[
\hatsetd_\rho(\epi f, \epi g^\lambda) \leq \big(1+\kappa(\hat\rho)\big)\max\Big\{\frac{\rho^*}{\lambda}, \psi^{-1}\Big(\frac{\rho^*}{\lambda}\Big)\Big\} + (1+m\lambda)\max_{i=0, \dots, m}\nsup_{\ball_X(\bar\rho)}|f_i-g_i|.
\]
as long as $\bar \rho > 2\rho + \max\{\dist((x^{\rm ctr},0), \epi f), \dist((x^{\rm ctr},0), \epi g^\lambda)\}$, $\rho^* \geq \bar \rho$ $+$ $\max\{0,$ $-\inf_{\ball_X(\bar\rho)} f_0\}$,  $\hat\rho > \bar\rho + \max\{\rho^*/\lambda, \psi^{-1}(\rho^*/\lambda)\}$, and the following constraint qualification holds: there is a strictly increasing function $\psi:\preals\to \preals$ such that
 \[
  \max_{i=1, \dots,m } f_i(x) \geq \psi\Big(\dist\big(x,\nlev_0 \{\max_{i=1, \dots, m} f_i\}\big)\Big) \mbox{ when } x\not\in \nlev_0 \{\max_{i=1, \dots, m} f_i\}.
  \]
\end{theorem}
\state Proof. As intermediate steps, we define $h,h^\lambda,f^\lambda:X\times\reals^m\to \Reals$ to have values $h(x,y) = \iota_{X\times \{0\}}(x,y) + \iota_C(x,y)$, with $C = \{(x,y)\in X\times \reals^m~|~f_i(x) \leq y_i, y_i\geq 0, i=1, \dots, m\}$, and
\[
h^\lambda(x,y) = \lambda \sum_{i=1}^m y_i + \iota_C(x,y)~~~~~f^\lambda(x,y) = f_0(x) + h^\lambda(x,y).
\]
First, we examine the Kenmochi conditions for $h$ and $h^\lambda$.
Let $(x,y) \in \nlev_{\rho^*} h^\lambda \cap \ball_{X\times \reals^m}(\rho^*)$. Thus, $(x,y)\in C$, $\lambda \sum_{i=1}^m y_i\leq \rho^*$, and $\|y\|_\infty \leq \rho^*/\lambda$. Let $\epsilon>0$ and $\eta = \max\{\rho^*/\lambda,\psi^{-1}(\rho^*/\lambda)\} + \epsilon$. If $f_i(x) \leq 0$ for all $i$, then
\[
\inf_{\ball_{X\times\reals^m}((x,y),\eta)} h \leq h(x,0) = 0 \leq \max\big\{h^\lambda(x,y), -\rho^*\big\}.
\]
Otherwise there is $i^*$ with $f_{i^*}(x) > 0$ so that
\[
\rho^*/\lambda\geq y_{i^*} \geq f_{i^*}(x) \geq \psi(\dist(x,\nlev_0\{ \max_{i=1, \dots, m} f_i\}))
\]
and $\psi^{-1}(\rho^*/\lambda) \geq \dist(x,\nlev_0 \{\max_{i=1, \dots, m} f_i\})$.  There exists $\bar x\in \nlev_0 \{\max_{i=1, \dots, m} f_i\}$ such that $d_X(x,\bar x) \leq \dist(x,\nlev_0 \{\max_{i=1, \dots, m} f_i\}) + \epsilon \leq \psi^{-1}(\rho^*/\lambda) \ + \epsilon$. Consequently,
\[
\inf_{\ball_{X\times\reals^m}((x,y),\eta)} h \leq h(\bar x,0) = 0 \leq \max\big\{h^\lambda(x,y), -\rho^*\big\}.
\]
Thus, the second set of Kenmochi conditions holds with this $\eta$.
Since $h^\lambda \leq h$, the first set also holds. Consequently, since $\epsilon>0$ is arbitrary and Prop. \ref{prop:altform} applies, we have establish that
\[
\hatsetd_{\rho^*}(\epi h, \epi h^\lambda) \leq \max\big\{\rho^*/\lambda,\psi^{-1}(\rho^*/\lambda)\big\}.
\]
We obtain via Prop. \ref{prop:lip} that
\[
\hatsetd_{\bar \rho}(\epi f, \epi f^\lambda) \leq \big(1+\kappa(\hat\rho)\big)\max\big\{\rho^*/\lambda, \psi^{-1}(\rho^*/\lambda)\big\}.
\]

Second, we consider the Kenmochi conditions for $f^\lambda$ and $g^\lambda$. Let $\delta = \max_{i=0, 1, \dots, m}$ $\sup_{\ball_X(\bar\rho)}|f_i-g_i|$ and $(x,y) \in \nlev_{\bar\rho} f^\lambda\cap\ball_{X\times\reals^m}(\bar\rho)$. Then, $(x,y)\in C$, $f_i(x) \leq y_i$, and $g_i(x) \leq y_i + \delta$ for all $i=1, \dots, m$. Set $\eta = (1+m\lambda)\delta$ and $\bar y = y + (\delta, \dots, \delta)$. With $B = \ball_{X\times\reals^m}((x,y),\eta)$, we obtain
\[
\ninf_{B} g^\lambda \leq g^\lambda(x, \bar y) = g_0(x) + \lambda\sum_{i=1}^m \bar y_i \leq f_0(x) + \delta + \lambda\sum_{i=1}^m y_i + \lambda m \delta \leq f^\lambda(x,y) + \eta.
\]
Repeating this argument with the roles of $g^\lambda$ and $f^\lambda$ reversed, we obtain via Prop. \ref{prop:altform} that $\hatsetd_{\bar\rho}(\epi f^\lambda, \epi g^\lambda) \leq (1+m\lambda)\delta$. Prop. \ref{prop:triangle} then yields the conclusion.\eop

The theorem presents a tradeoff between two error terms. If the conditioning function $\psi(\gamma) = \gamma^\beta$ for $\beta>0$, then $\lambda$ should be of the order $O(\delta^{-\beta/(1+\beta)})$ to balance the two terms, where $\delta = \nmax_{i=0, 1, \dots, m}\nsup_{\ball_X(\bar\rho)}|f_i-g_i|$. This leads to the overall rate of convergence $O(\delta^{1/(1+\beta)})$, which can be significantly worse than what is indicated by the pointwise error $\delta$. Still, the situation is much improved from the approach of simply minimizing $g_0(x)$ subject to $g_i(x) \leq 0$ for $i=1, \dots, m$. As discussed prior to the theorem, that problem may have solutions that are arbitrarily far away from those of the actual problem \eqref{eqn:actualproblem}. In some sense, the theorem explains the popularity of formulations with constraint softening in practice (see \cite{BrownCarlyle.08} for a prime example); they are in a fundamental way ``robust'' to inaccuracy in the constraint functions.

Theorem \ref{thm:slack} makes no Slater-type constraint qualification for the actual problem and places no restrictions on the properties of the constraint functions at points in the feasible set. Naturally, if such conditions are brought in, we can improve the results; cf. Prop. \ref{prop:convexlevel} and \cite[Thm. 4.6]{Royset.18}.\\

\noindent {\bf Case III.} While still addressing the actual problem \eqref{eqn:actualproblem}, the third case examines the classical penalty method and the resulting unconstrained approximating problems.

\begin{theorem}\label{thm:penality}{\rm (approximation by penalty formulation).} For a metric space $X$, with centroid $x^{\rm ctr}$, $\lambda\in (0, \infty)$, and $f_i, g_i:X\to \reals$, $i=0, 1, \dots, m$, where $f_0$ and $g_0$ are Lipschitz continuous with common modulus $\kappa:\preals\to \preals$, consider the functions $f,g^\lambda:X\times \reals^m\to \Reals$ defined by
\[
g^\lambda(x) = g_0(x) + \lambda \sum_{i=1}^m \max\{0, g_i(x)\} \mbox{ and } f(x) = \begin{cases}
f_0(x) &\mbox{ if } f_i(x) \leq 0 ~\forall i = 1, \dots, m\\
\infty &\mbox{ otherwise}
\end{cases}
\]
Then,
\[
\hatsetd_\rho(\epi f, \epi g^\lambda) \leq \max\{1,\kappa(\hat\rho)\}\psi^{-1}\Big(\frac{\bar \rho - \ninf_{\ball_X(\bar\rho)} f_0}{\lambda}\Big) + (1+m\lambda)\max_{i=0, \dots, m}\sup_{\ball_X(\bar\rho)}|f_i-g_i|.
\]
provided that $\bar \rho > 2\rho + \max\{\dist(x^{\rm ctr}, \epi f), \dist(x^{\rm ctr}, \epi g^\lambda)\}$, $\hat\rho > \bar\rho + \psi^{-1}((\bar \rho - \ninf_{\ball_X(\bar\rho)} f_0)\lambda^{-1})$, and the same constraint qualification as in Thm. \ref{thm:slack} holds.
\end{theorem}
\state Proof. As an intermediate quantity, we define $f^\lambda:X\to\reals$ to have values $f^\lambda(x) = f_0(x) + \lambda \sum_{i=1}^m \max\{0, f_i(x)\}$. We start by examining the Kenmochi conditions for  $f$ and $f^\lambda$. Let $x \in \nlev_{\bar \rho} f^\lambda \cap \ball_X(\bar\rho)$ so that $f_0(x) + \lambda \sum_{i=1}^m \max\{0, f_i(x)\} \leq \bar\rho$. If $\max_{i=1, \dots, m} f_i(x) >0 $, then
\[
\nmax_{i=1, \dots, m} f_i(x) \leq \sum_{i=1}^m \max\{0, f_i(x)\} \leq \frac{\bar\rho -f_0(x)}{\lambda}.
\]
Since $f_0(x) \leq \bar\rho$, $\inf_{\ball_X(\bar\rho)} f_0 \leq \bar\rho$. These facts together with the constraint qualification lead to
\[
\dist\big(x, \nlev_0 \{\max_{i=1, \dots, m} f_i\}\big) \leq \psi^{-1}\big(\max_{i=1, \dots, m} f_i(x)\big) \leq \eta = \psi^{-1}\left(\frac{\bar\rho - \inf_{\ball_X(\bar\rho)} f_0}{\lambda}\right).
\]
Let $\epsilon\in(0, \hat\rho - \bar\rho - \psi^{-1}((\bar \rho - \ninf_{\ball_X(\bar\rho)} f_0)\lambda^{-1})]$. There exists $\bar x\in \nlev_0 \{\nmax_{i=1, \dots, m} f_i\}$ such that $d_X(x,\bar x) \leq \eta + \epsilon$ and
\[
\ninf_{\ball_X(x,\eta+\epsilon)} f \leq f(\bar x) = f_0(\bar x) \leq f_0(x) + \kappa(\hat\rho)(\eta+\epsilon) \leq f^\lambda (x)+ \kappa(\hat\rho)(\eta+\epsilon).
\]
Alternatively, if $\max_{i=1, \dots, m} f_i(x) \leq 0$, then $\ninf_{\ball_X(x,0)} f \leq f_0(x) \leq f^\lambda (x)$. We have therefore established the second Kenmochi condition for  $f$ and $f^\lambda$ with error $\max\{1, \kappa(\hat\rho)\}(\eta+ \epsilon)$. Since $f\geq f^\lambda$, the first Kenmochi condition holds with an error of zero. Since $\epsilon>0$ is arbitrary, we have established via Prop. \ref{prop:altform} that
\[
\hatsetd_{\bar\rho}(\epi f,\epi f^\lambda) \leq \max\big\{1,\kappa(\hat\rho)\big\} \psi^{-1}\left(\frac{\bar\rho - \inf_{\ball_X(\bar\rho)} f_0}{\lambda}\right).
\]
Trivially, $|f^\lambda(x) - g^\lambda(x)| \leq (1+m\lambda)\nmax_{i=0, 1, \dots, m}\nsup_{\ball_X(\bar\rho)}|f_i-g_i|$ for $x\in \ball_X(\bar\rho)$ so that  $\hatsetd_{\bar\rho}(\epi f^\lambda, \epi g)$ is also bounded by the same quantity; cf. Prop. \ref{prop:fromsupnorm}. The conclusion then follows by Prop. \ref{prop:triangle}.\eop

We again find a tradeoff between two error terms that are nearly identical to those in Thm. \ref{thm:slack}. From this perspective, the penalty formulation has the same rate of convergence as that in Case II and is therefore stable even when the actual feasible set in \eqref{eqn:actualproblem} has an empty interior.

\subsection{Calculus Rules for Compositions}

The truncated Hausdorff distance between epigraphs of functions that are certain compositions can be bounded as we see next. The results of this subsection extend in some sense Prop. \ref{prop:lip}, which deals with sums. Composition rules for epi-sum and epi-multiplication can be found in \cite{AttouchWets.91}; see also \cite{AzePenot.90b} for a systematic treatment of the convex case including sums of convex functions.

\begin{proposition}\label{prop:compLipsch}{\rm (compositions; Lipschitz inner mapping).} For metric spaces $(X,d_X)$ and $(Y, d_Y)$, with centroids $x^{\rm ctr}$ and $y^{\rm ctr}$, respectively, $f,g:Y\to \Reals$, and $F,G:X\to Y$, suppose that $F^{-1},G^{-1}:Y\tto X$ are nonempty-valued and Lipschitz continuous with common modulus $\kappa:\preals\to \preals$ relative to $\rho^*\in [0,\infty]$. Then, for $\rho\in \preals$,
\[
\hatsetd_\rho\big( \epi (f \circ F), \epi (g \circ G)\big) \leq \sup_{y\in \ball_{Y}(\bar\rho)} \hatsetd_{\rho^*}(F^{-1}(y),G^{-1}(y))    + \max\{1,\kappa(\hat \rho)\} \hatsetd_{\bar \rho}(\epi f, \epi g)
\]
provided that $\rho^* > 2\rho + \max\{|\alpha|, |\bar\alpha|, \dist(x^{\rm ctr}, F^{-1}(y)), \dist(x^{\rm ctr}, F^{-1}(\bar y)), and
\newline
\dist(x^{\rm ctr}, G^{-1}(\bar y))\}$
for some $(y,\alpha)\in \epi f$ and $(\bar y, \bar\alpha) \in \epi g$,
\[
\bar \rho >  \max\big\{\rho^*, ~\nsup_{x\in\ball_X(\rho^*)} d_Y\big(F(x) ,y^{\rm ctr}\big), ~\nsup_{x\in\ball_X(\rho^*)} d_Y\big(G(x) ,y^{\rm ctr}\big)\big\},
\]
and $\hat \rho > \bar \rho + \hatsetd_{\bar\rho}(\epi f, \epi g)$.
\end{proposition}
\state Proof. Let $\hat F, \hat G:X\times\reals \to Y\times \reals$ have $\hat F(x,\alpha) = (F(x), \alpha)$ and $\hat G(x,\alpha) = (G(x), \alpha)$ for $(x,\alpha)\in X\times \reals$. Then, it follows directly that
\[
\epi (f\circ F) = \hat F^{-1}(\epi f)~ \mbox{ and } ~\epi (g\circ G) = \hat G^{-1}(\epi g)
\]
and we can bring in Thm. \ref{thm:Lipmapping} with $S = \hat F^{-1}$ and $T = \hat G^{-1}$. Let $\epsilon> 0$. There exists $x\in F^{-1}(y)$ such that $d_X(x^{\rm ctr}, x) \leq \dist(x^{\rm ctr}, F^{-1}(y)) + \epsilon$. Then, $f(F(x)) = f(y) \leq \alpha$ and $(x,\alpha) \in \epi (f\circ F)$. Consequently,
\[
\dist\big((x^{\rm ctr}, 0), \epi (f\circ F) \big) \leq  \max\big\{d_X(x^{\rm ctr}, x), |\alpha|\big\} \leq \max\big\{\dist(x^{\rm ctr}, F^{-1}(y)) + \epsilon, |\alpha|\big\}.
\]
Similar arguments establish that
\begin{align*}
\dist\big((x^{\rm ctr}, 0), \epi (g\circ F) \big) &\leq  \max\big\{\dist(x^{\rm ctr}, F^{-1}(\bar y)) + \epsilon, |\bar\alpha|\big\}\\
\dist\big((x^{\rm ctr}, 0), \epi (g\circ G) \big) &\leq  \max\big\{\dist(x^{\rm ctr}, G^{-1}(\bar y)) + \epsilon, |\bar\alpha|\big\}.
\end{align*}
This ensures that $\rho^*$ is selected sufficiently large for the application of Thm. \ref{thm:Lipmapping}. Next, we consider the size of $\bar \rho$ and find that
\begin{align*}
&\sup_{(x,\alpha)\in \hat F^{-1}(\epi f) \cap\ball_{X\times \reals}(\rho^*)} \Big\{\inf_{\hat F(x,\alpha)\cap \epi f} d_{Y\times\reals} \big(\cdot, (y^{\rm ctr}, 0)\big)\Big\}\\
= &\sup\Big\{ \max\big\{d_Y(F(x),y^{\rm ctr}), |\alpha|\big\}   ~\Big|~f\big(F(x)\big)\leq \alpha, x\in \ball_X(\rho^*), |\alpha|\leq \rho^*   \Big\}\\
\leq &\max\Big\{\rho^*, \nsup_{x\in \ball_X(\rho^*)} d_Y\big(F(x),y^{\rm ctr}\big)\Big\}.
\end{align*}
Since similar statements hold with $F$ replaced by $G$ and $\epi f$ replaced by $\epi g$, the condition on $\bar \rho$ suffices and Thm. \ref{thm:Lipmapping} yields the conclusion.\eop

\begin{corollary}\label{cor:compAffine}{\rm (compositions; linear inner mapping).} For $f,g:\reals^n\to \Reals$ and nonsingular $n\times n$ matrices $A$ and $B$, suppose that $\phi,\psi:\reals^n\to \Reals$ are defined by $\phi(x) = f(Ax)$ and $\psi(x) = g(Bx)$, $x\in \reals^n$. Then\footnote{Here we use the operator norm for matrices.}, for $\rho\in \preals$,
\[
\hatsetd_\rho( \epi \phi, \epi \psi) \leq \bar\rho\|A^{-1}-B^{-1}\| + \max\big\{1,\|A^{-1}\|, \|B^{-1}\|\big\} \hatsetd_{\bar \rho}(\epi f, \epi g)
\]
as long as $\bar \rho >  \max\{1, \|A\|, \|B\|\}(2\rho + \max\{|\alpha|, |\bar\alpha|, \dist(0, A^{-1}y)$, $\dist(0, A^{-1}\bar y)$, and $\dist(0, B^{-1}\bar y)\})$ for some $(y,\alpha)\in \epi f$ and $(\bar y, \bar\alpha) \in \epi g$.
\end{corollary}
\state Proof. The result follows directly from Prop. \ref{prop:compLipsch}.\eop

The corollary extends in some sense \cite[Cor. 2.6]{AzePenot.90b} by allowing for nonconvex $f,g$ and different linear mappings, but at the expense of requiring invertible mappings.

\begin{proposition}{\rm (compositions; Lipschitz outer function).} For metric spaces $(X,d_X)$ and $(Y, d_Y)$, with $y^{\rm ctr}$ being the centroid of $Y$, suppose that $f:Y\to \reals$ is Lipschitz continuous with modulus $\kappa:\preals\to \preals$, and $F,G:X\to Y$. Then, for $\rho\in \preals$,
\[
\hatsetd_\rho\big( \epi (f \circ F), ~\epi (f \circ G)\big) \leq \max\big\{1, \kappa(\hat\rho)\big\}\hatsetd_{\bar\rho}(\gph F, \gph G)
\]
provided that $\hat\rho > \bar\rho + \hatsetd_{\bar\rho}(\gph F, \gph G)$ and
\begin{align*}
\bar \rho > \max\Big\{\rho, &\sup_{x\in \ball_X(\rho)}\big\{d_Y(y^{\rm ctr}, F(x))~\big|~f(F(x))\leq \rho\big\},\\
&\sup_{x\in \ball_X(\rho)}\big\{d_Y(y^{\rm ctr}, G(x))~\big|~f(G(x))\leq \rho\big\} \Big\}.
\end{align*}
\end{proposition}
\state Proof. Let $\eta = \hatsetd_{\bar\rho}(\gph F, \gph G)$, $x\in \nlev_\rho (f\circ F) \cap \ball_X(\rho)$, and $\epsilon\in (0, \hat\rho - \bar\rho - \eta]$.  Then, $(x,F(x))\in \ball_{X\times Y}(\bar \rho)$ and there exists $\bar x\in X$ with $d_X(\bar x, x) \leq \eta + \epsilon$ and $d_Y(F(x),G(\bar x))\leq \eta + \epsilon$. Since both $F(x),G(\bar x) \in \ball_Y(\hat \rho)$,
\[
\ninf_{\ball_X(x,\eta + \epsilon)} (f\circ G) \leq f\big(G(\bar x)\big) \leq f\big(F(x)\big) + \kappa(\hat\rho)(\eta+ \epsilon).
\]
We repeat the argument with the roles of $F$ and $G$ reversed and obtain via Prop. \ref{prop:altform} that $\hatsetd_\rho(\epi (f\circ F), \epi (f\circ G)) \leq \max\{1, \kappa(\hat\rho)\}(\eta + \epsilon)$. Since $\epsilon$ is arbitrary, the conclusion follows.\eop

The previous two propositions largely summarize  the line of reasoning in the proofs of Thm.  \ref{thm:feasproblem}, \ref{thm:slack}, and \ref{thm:penality} and thereby facilitate various extensions of Cases I, II, and III.

\begin{proposition}\label{prop:infprojection}{\rm (inf-projections).}
For a metric space $X$ and $\{f_\alpha, g_\alpha:X\to \Reals\}$, with $A$ an arbitrary set, define $f,g:\reals^n\to \Reals$ as $f(x) = \ninf_{\alpha\in A} f_\alpha(x)$ and $g(x) = \ninf_{\alpha\in A} g_\alpha(x)$. Then, for $\rho\in \preals$,
\[
\hatsetd_\rho(\epi f, \epi g) \leq \nsup_{\alpha\in A} \hatsetd_\rho(\epi f_\alpha, \epi g_\alpha).
\]
\end{proposition}
\state Proof.  In view of the fact that $\epi f = \cup_{\alpha \in A} \epi f_\alpha$ and similarly for $\epi g$, the conclusion follows immediately from Prop. \ref{prop:union}.\eop

Since a function $f = \sup_{\alpha\in A} f_\alpha$ has as epigraph the intersection of $\epi f_\alpha, \alpha\in A$, it is clear from the discussion in Section 3 that no comparable result is possible for sup-projections. We refer to \cite[Cor. 2.5]{AzePenot.90b} for a result in the convex case and \cite[Thm. 5.6]{RoysetWets.16b} for one under Lipschitz continuity assumptions.

Given metric spaces $X$ and $Y$ as well as $f:X\to \Reals$ and $F:X\to Y$, the {\it epi-composition} $Ff:Y\to \Reals$ has
\[
(Ff)(y) := \ninf\big\{f(x)~|~F(x) = y\big\} \mbox{ for  } y\in Y.
\]
Epi-compositions arise, for example, in parametric studies of equality constrained problems.

\begin{proposition}\label{prop:epicomp}{\rm (epi-compositions).}
For metric spaces $(X,d_X)$ and $(Y,d_Y)$, with $x^{\rm ctr}$ being the centroid of $X$, $f,g:X\to\reals$, and Lipschitz continuous $F,G:X\to Y$ with common modulus $\kappa:\preals \to \preals$ relative to $\infty$, suppose that
\begin{align*}
y\in Y \mbox{ and } (Ff)(y) \in \reals &\mbox{ imply } \nargmin_{x\in X} \{f(x)~|~F(x) = y\} \neq\emptyset; \mbox{ and }\\
y\in Y \mbox{ and } (Gg)(y) \in \reals &\mbox{ imply } \nargmin_{x\in X} \{g(x)~|~G(x) = y\} \neq\emptyset.
\end{align*}
Then, for $\rho\in \preals$,
\[
\hatsetd_\rho(\epi Ff, \epi Gg) \leq \nsup_{x\in \ball_X(\bar \rho)}d_Y\big(F(x),G(x)\big) + \max\{1, \kappa(\hat \rho)\}\hatsetd_{\bar\rho}(\epi f, \epi g)
\]
provided that $\rho^* > 2\rho + \max\{d_Y(F(x),y^{\rm ctr}),$ $d_Y(F(\bar x),y^{\rm ctr}), d_Y(G(\bar x),y^{\rm ctr}), |\alpha|, |\bar \alpha|\}$ for some $(x,\alpha) \in \epi f$ and $(\bar x, \bar \alpha) \in \epi g$, $\bar\rho>\rho^*$ and also exceeds
\[
\sup\Big\{d_X(x,x^{\rm ctr})~\Big| (x,\alpha)\in C, |\alpha| \leq \rho^*, U(x) \in \ball_Y(\rho^*)\Big\} \mbox{ for } U = F,G; C = \epi f, \epi g,
\]
and $\hat \rho > \bar \rho + \hatsetd_{\bar\rho}(\epi f, \epi g)$.
\end{proposition}
\state Proof. We start by confirming that $\epi Ff = \{(F(x),\alpha)~|~(x,\alpha) \in \epi f\}$; a finite-dimensional version of this fact is asserted as Exercise 1.31 in \cite{VaAn}. For $(\bar x, \bar\alpha)\in \epi f$, we have that $\inf\{f(x)~|~F(x) = F(\bar x)\} \leq f(\bar x) \leq \bar\alpha$. Thus, $\epi Ff \supset \{(F(x),\alpha)~|~(x,\alpha)$ $\in$ $\epi f\}$. Suppose that $(y,\alpha) \in \epi Ff$. Then, $(Ff)(y) < \infty$. If $(Ff)(y)= -\infty$, then there exists $\bar x\in X$ such that $f(\bar x) \leq \alpha$ and $F(\bar x) = y$. Consequently, $(y,\alpha) \in  \{(F(x),\alpha)~|~(x,\alpha) \in \epi f\}$. If $(Ff)(y) \in \reals$, then there exists by assumption $\bar x\in X$ such that $f(\bar x) = \inf\{f(x)~|~F(x) = y\}$ and $F(\bar x) = y$. Thus, $f(\bar x) = (Ff)(y) \leq \alpha$, $(\bar x, \alpha) \in \epi f$, and $\epi Ff \subset \{(F(x),\alpha)~|~(x,\alpha) \in \epi f\}$. We have confirmed the assertion, which also holds for $Gf$.

The conclusion follows by Thm. \ref{thm:Lipmapping} applied to the mappings $\hat F,\hat G:X\times\reals\to Y\times \reals$ defined by $\hat F(x,\alpha) = (F(x),\alpha)$ and $\hat G(x,\alpha)=(G(x),\alpha)$. Since $F$ and $G$ are Lipschitz continuous with common modulus $\kappa:\preals\to\preals$ relative to $\infty$, $\hat F$ and $\hat G$ are Lipschitz continuous with modulus $\rho\mapsto \max\{1, \kappa(\rho)\}$ relative to any real number. The requirement on $\rho^*$ in Thm. \ref{thm:Lipmapping} is satisfied because $\dist((y^{\rm ctr}, 0), \hat F(\epi f)) \leq \max\{d_Y(F(x), y^{\rm ctr}), |\alpha|\}$ for $(x,\alpha) \in \epi f$, with similar inequalities holding for $\hat G$ and $\epi g$. The requirement on $\bar\rho$ in Thm. \ref{thm:Lipmapping} also is satisfied because
\begin{align*}
&\sup_{(y,\alpha)\in \hat F(\epi f) \cap \ball_{Y\times\reals}(\rho^*)} \Big\{\inf_{\hat F^{-1}(y,\alpha) \cap \epi f} d_X\big(\cdot,(x^{\rm ctr},0)\big)\Big\}\\
 \leq &\sup_{(x,\alpha)\in \epi f, |\alpha| \leq \rho^*, F(x) \in \ball_Y(\rho^*)} \max\{d_X(x,x^{\rm ctr}), |\alpha|\}
\end{align*}
with similar expressions for $\hat G$ and $\epi g$.\eop

\section{Distances between Graphs of Set-Valued Mappings}

We next turn to the solution of generalized equations. For metric spaces $X$ and $Y$, a set-valued mapping $S:X\tto Y$ and a point $y^\star\in Y$ define the generalized equation $y^\star \in S(x)$. Its {\it solution set} is $S^{-1}(y^\star)$. In this section, we focus on the  {\it set of near-solutions} that consists of those $x\in X$ with $S(x)$ ``nearly reaching'' $y^\star$. Specifically, for $\epsilon\geq 0$, the {\it set of $\epsilon$-solutions} is defined as
\[
S^{-1}\big( \ball_Y(y^\star, \epsilon)  \big) = \bigcup_{y\in \ball_Y(y^\star, \epsilon)} S^{-1}(y).
\]

For example, suppose that $f:\reals^n\to \reals$ is locally Lipschitz continuous and $C\subset \reals^n$ is nonempty and closed. Then, an optimality conditions for the problem of minimizing $f + \iota_C$ would be
\[
0 \in \partial f(x) + N_C(x);
\]
see \cite[Exercise 10.10]{VaAn}. With $S = \partial f + N_C$ and $y^\star=0$, the set of $\epsilon$-solutions becomes
\[
S^{-1}\big( \ball_{\reals^n}(\epsilon)  \big) = \big\{x\in \reals^n~|~ 0 \in \partial f(x) + N_C(x) + \ball_{\reals^n}(\epsilon)\big\}.
\]
The next theorem bounds the discrepancy between near-solutions of generalized equations in terms of the truncated Hausdorff distance without making assumptions about local regularity properties of the underlying set-valued mappings.

\begin{theorem}\label{thm:approxgeneralizedequations}{\rm (approximation of near-solutions of generalized equations).} For metric spaces $X$ and $Y$, suppose that $S,T:X\tto Y$ have nonempty graphs, $0 \leq\epsilon\leq \rho< \infty$, and $y^\star\in \ball_Y(\rho-\epsilon)$. Then,
\[
\exs\Big(S^{-1}\big(\ball_Y(y^\star, \epsilon)\big)\cap \ball_X(\rho); ~T^{-1}\big(\ball_Y(y^\star, \delta)\big)    \Big) \leq \hatsetd_\rho(\gph S, \gph T)
\]
provided that $\delta > \epsilon + \hatsetd_\rho(\gph S, \gph T)$. If $X$ and $Y$ are finitely compact and $\gph T$ is closed, then the result also holds for $\delta = \epsilon + \hatsetd_\rho(\gph S, \gph T)$.
\end{theorem}
\state Proof. Let $\gamma\in (0, \delta - \epsilon - \hatsetd_\rho(\gph S, \gph T)]$. Suppose that $x\in S^{-1}\big(\ball_Y(y^\star, \epsilon)\big)\cap \ball_X(\rho)$. Then, there is $y\in S(x)$ with $d_Y(y,y^\star)\leq \epsilon$ so that $(x,y) \in \ball_{X\times Y}(\rho)$.  Consequently, for some $(\bar x, \bar y)\in \gph T$,
\[
 \max\big\{d_X(x, \bar x), d_Y(y, \bar y) \big\} \leq  \dist\big((x, y), \gph T\big) + \gamma \leq \hatsetd_\rho(\gph S, \gph T)+\gamma.
\]
Moreover, $d_Y(\bar y, y^\star) \leq     d_Y(\bar y, y) + d_Y(y, y^\star) \leq \hatsetd_\rho(\gph S, \gph T)+\gamma + \epsilon \leq \delta$, which implies that $\bar x \in  T^{-1}(\ball_Y(y^\star, \delta))$. We have established that
\[
\exs\Big(S^{-1}\big(\ball_Y(y^\star, \epsilon)\big)\cap \ball_X(\rho); ~T^{-1}\big(\ball_Y(y^\star, \delta)\big)    \Big) \leq \hatsetd_\rho(\gph S, \gph T) + \gamma.
\]
Since $\gamma$ is arbitrary, the first conclusion follows. The minimum distance to a nonempty closed subset of a finitely compact space is attained \cite[Lemma 2.2]{Royset.18}, which allows us to use $\gamma = 0$ in the above arguments. This establishes the second conclusion.\eop

The result of the theorem is sharp. For example, consider $S,T:\reals\tto\reals$ with $S(x) = [x,\infty)$ when $x\in [0,1]$ and $S(x) = \emptyset$ otherwise; and $T(x) = (1,\infty)$ when $x\in [1,2]$ and $T(x) = \emptyset$ otherwise. Then for $\rho \geq 0$, $\hatsetd_\rho(\gph S, \gph T) = 1$, $S^{-1}(0) = \{0\}$, $T^{-1}(\delta)=[1,2]$, and $\exs(S^{-1}(0)\cap \ball_{\reals}(\rho); T^{-1}(\ball_\reals(\delta)) = 1$  when $\delta > 1$. When $\delta \leq 1$, the excess becomes infinity because $T^{-1}(\delta) = \emptyset$. If $T$ is modified to having $T(x) = [1,\infty)$ for $x\in [1,2]$, then $\delta = 1$ gives an excess of one.

\begin{theorem}\label{thm:summappings}{\rm (sum of mappings under Lipschitz property).} For normed linear spaces $X$ and $Y$, suppose that $S_1, T_1:X\tto Y$ are nonempty-valued and Lipschitz continuous with common modulus $\kappa:\preals\to \preals$ relative to $\rho^* \in [0,\infty]$ and $S_2, T_2:X\tto Y$ have nonempty graphs. Then, for $\rho\in \preals$,
\begin{align*}
\hatsetd_\rho\big(\gph (S_1 + S_2), \gph (T_1 + T_2)\big) \leq &\nsup_{x\in \ball_{X}(\rho)} \hatsetd_{\rho^*}\big(S_1(x), T_1(x)\big)\\
& + \big(1 + \kappa(\hat\rho)\big)\hatsetd_{\bar \rho}(\gph S_2, \gph T_2),
\end{align*}
provided that $\bar \rho \geq  \rho + \rho'$, with $\rho'$ such that $\ball_{Y}(\rho')$ contains both $S_1(x)$ and $T_1(x)$ for all $x\in \ball_{X}(\rho)$,  $\hat\rho > \rho + \hatsetd_{\bar \rho}(\gph S_2, \gph T_2)$, and $\rho^* > 3\rho' + \kappa(\hat\rho)(\hat\rho-\rho)$.
\end{theorem}
\state Proof. Let $(x,y) \in \gph (T_1 + T_2) \cap \ball_{X \times Y}(\rho)$. Thus, for some $y_1\in T_1(x)$ and $y_2\in T_2(x)$ we have $y = y_1 + y_2$ and $\|y_2\| \leq \|y\| + \|y_1\| \leq \rho + \rho'\leq \bar\rho$. Let $\epsilon\in (0, \hat\rho - \rho -\hatsetd_{\bar \rho}(\gph S_2, \gph T_2)]$. Consequently, $(x,y_2) \in \gph T_2 \cap \ball_{X\times Y}(\bar \rho)$ so there exists $(\bar x, \bar y_2)\in \gph S_2$ with $\max\{\|x - \bar x\|, \|y_2 - \bar y_2\|\}\leq \hatsetd_{\bar \rho}(\gph S_2, \gph T_2) + \epsilon \leq \hat\rho - \rho$, which ensures that $\|\bar x\| \leq \|x-\bar x\| + \|x\| \leq \hat\rho - \rho + \rho \leq \hat\rho$. Since $S_1$ is nonempty-valued, there is $\bar y_1 \in S_1(\bar x)$  such that $\dist(y_1, S_1(\bar x)) \geq \|y_1 - \bar y_1\|-\epsilon$.
Therefore,  $(\bar x, \bar y_1 + \bar y_2) \in \gph (S_1 + S_2)$.
Since $y_1\in \ball_Y(\rho')$, it follows  that
\begin{align*}
\|y_1 - \bar y_1\| &\leq  \dist\big(y_1, S_1(\bar x)\big) + \epsilon \leq \hatsetd_{\rho'}\big(S_1(\bar x), T_1(x)\big)  +\epsilon\\
&\leq \hatsetd_{\rho^*}\big(S_1(\bar x), S_1(x)\big)  + \hatsetd_{\rho^*}\big(S_1(x), T_1(x)\big)  + \epsilon,
\end{align*}
where the last inequality is a consequence of Prop. \ref{prop:triangle}; $\rho^*$ is indeed sufficiently large because $\dist(y^{\rm ctr}, T_1(x))\leq \rho'$, $\dist(y^{\rm ctr}, S_1(x))\leq \rho'$, and
\begin{align*}
\dist\big(y^{\rm ctr}, S_1(\bar x)\big)&\leq \rho' + \exs\big(S_1(x) \cap \ball_Y(\rho^*); S_1(\bar x)\big) \leq \rho' + \hatsetd_{\rho^*}\big(S_1(x), S_1(\bar x)\big)\\
& \leq \rho' + \kappa(\hat\rho)\|x-\bar x\| \leq \rho' + \kappa(\hat\rho)(\hat\rho - \rho).
\end{align*}
Moreover, with $\bar y = \bar y_1 + \bar y_2$, $\|y - \bar y\|$ is not greater than
\begin{align*}
  &\|y_1 - \bar y_1\| + \|y_2 - \bar y_2\| \leq \hatsetd_{\rho^*}\big(S_1(\bar x), S_1(x)\big)\\
   & ~~~~~~~~~~~~~~~~~~~~+ \hatsetd_{\rho^*}\big(S_1(x), T_1(x)\big)  + \hatsetd_{\bar \rho}(\gph S_2, \gph T_2) + 2\epsilon\\
   &\leq \kappa(\hat\rho)\big[\hatsetd_{\bar \rho}(\gph S_2, \gph T_2) + \epsilon\big] + \nsup_{x' \in \ball_{X}(\rho)} \hatsetd_{\rho^*}\big(S_1(x'), T_1(x')\big)\\
   & ~~~~~~~~~~~~~~~~~~~~+ \hatsetd_{\bar \rho}(\gph S_2, \gph T_2) + 2\epsilon.
\end{align*}
This establishes that $(\bar x, \bar y) \in \gph (S_1+S_2)$ satisfies
\begin{align*}
&\max\{\|x - \bar x\|, \|y-\bar y\|\} \leq \max\big\{\hatsetd_{\bar \rho}(\gph S_2, \gph T_2) + \epsilon,\\
&~~~~~~~~~ \big(1 + \kappa(\hat\rho)\big)\hatsetd_{\bar \rho}(\gph S_2, \gph T_2)   + \nsup_{x'
\in\ball_{X}(\rho)} \hatsetd_{\rho^*}\big(S_1(x'), T_1(x')\big)  + (2 + \kappa(\hat\rho))\epsilon\big\}.
\end{align*}
Since $(x,y)$ and $\epsilon$ are arbitrary, we obtain that
\begin{align*}
&\exs\big(\gph (T_1 + T_2) \cap \ball_{X \times Y}(\rho); \gph (S_1+S_2)  \big)\\
& \leq \nsup_{x' \in \ball_{X}(\rho)} \hatsetd_{\rho^*}\big(S_1(x'), T_1(x')\big) + \big(1 + \kappa(\hat\rho)\big)\hatsetd_{\bar \rho}(\gph S_2, \gph T_2).
\end{align*}
The roles of $(S_1, S_2)$ and $(T_1, T_2)$ can be reversed, which leads to the conclusion.\eop

A series of results are now possible with applications to games as well as equilibrium and generalized fixed-point problems. We limit the discussion to optimality conditions. As a preliminary example, let  $C,D\subset \reals^n$ be nonempty, possibly nonconvex sets and $f,g:\reals^n\to \reals$ be smooth and their gradients be Lipschitz continuous with modulus $\kappa:\preals\to \preals$ relative to $\rho^* = \infty$, i.e., $\|\nabla f(x) - \nabla f(\bar x)\| \leq \kappa(\rho)\|x - \bar x\|$ for $\|x\|\leq \rho$, $\|\bar x\|\leq \rho$, and $\rho \in \preals$, with the same condition holding for $\nabla g$. Thm. \ref{thm:summappings} enables a study of the optimality conditions $0 \in \nabla f(x) + N_C(x)$ and $0 \in \nabla g(x) + N_D(x)$. The discrepancy between the corresponding near-stationary points are bounded via Thm. \ref{thm:approxgeneralizedequations} by
\begin{align*}
\hatsetd_\rho\big(\gph(\nabla f + N_C), \gph(\nabla g + N_D)\big) \leq &\nsup_{\|x\|\leq\rho} \|\nabla f(x) - \nabla g(x)\|\\
& + \big(1+\kappa(\hat\rho)\big) \hatsetd_{\bar\rho } \big(\gph N_C, \gph N_D\big)
\end{align*}
for sufficiently large $\hat\rho$ and $\bar\rho$ with further simplifications possible if $C$ and $D$ are convex, cf. Prop. \ref{prop:normalcone}.\\

{\it Example 3: difference-of-convex functions}. For convex functions $f_1:\reals^n\to \reals$ and  $f_2:\reals^n\to \Reals$, the latter also lsc and proper,  as well as a point $\bar x$ with $f_2(\bar x)$ finite, the following optimality condition holds\footnote{For subsets $A$ and $B$ of a linear space, $A-B := \{a-b~|~a\in A, b\in B\}$.} \cite{HiriartUrruty.85}:
\[
\bar x \mbox{ local minimizer of } f_2 - f_1 \Longrightarrow 0 \in \partial f_2(\bar x) - \partial f_1(\bar x).
\]
The minimization of such difference-of-convex functions arises in numerous applications include some in moderns statistics \cite{CuiPangSen.18,Royset.19}. Error analysis of near-stationarity in this case can be carried our as follows.

Suppose initially that $f_1,g_1$ are also smooth and $\rho\in \preals$. Then, there are $\alpha,\bar\rho\in \preals$ such that\footnote{We here use the Euclidean distance on $\reals^n$.}
\begin{align*}
\hatsetd_\rho\big(\gph (\partial f_2 - \nabla f_1), \gph (\partial g_2 - \nabla g_1)\big) \leq & \nsup_{\|x\|_2\leq \rho} \|\nabla f_1(x) - \nabla g_1(x)\|_2\\
 + \alpha\sqrt{\hatsetd_{\bar\rho}(\epi f_2, \epi g_2 ) },
\end{align*}
which via Thm. \ref{thm:approxgeneralizedequations} gives error estimates of near-stationary points. We can establish this fact by setting $S_1 = -\nabla f_1$, $T_1 = -\nabla g_1$, $S_2 = \partial f_2$, and $T_2 = \partial g_2$ so that $S_1$ and $T_1$ are nonempty-valued and Lipschitz continuous with some common modulus $\kappa:\preals\to\preals$ relative to $\rho^* = \infty$. An application of Thm. \ref{thm:summappings} with these set-valued mappings and $\rho' = \sup_{\|x\|_2\leq \rho} \max\{\|\nabla f(x)\|_2,$ $\|\nabla g(x)\|_2\}$, $\bar\rho = \rho + \rho'$, and $\hat\rho > \rho + \hatsetd_{\bar \rho}(\gph \partial f_2, \gph \partial g_2)$ yields
\begin{align*}
\hatsetd_\rho\big(\gph (\partial f_2 - \nabla f_1), \gph (\partial g_2 - \nabla g_1)\big) \leq & \nsup_{\|x\|_2\leq \rho} \|\nabla f_1(x) - \nabla g_1(x)\|_2\\
& + \big(1 + \kappa(\hat\rho)\big)\hatsetd_{\bar \rho}(\gph \partial f_2, \gph \partial g_2 ).
\end{align*}
An application of Prop. \ref{prop:subgradmappings} gives the result after an appropriate enlargement of $\bar\rho$.

We can relax the assumption about $f_1$ and $g_1$ being smooth by stating the optimality condition in terms of the set-valued mappings $S,T:\reals^n\times\reals^n\tto \reals^n\times\reals^n$ with expressions
\[
S(x,v) = \begin{pmatrix}
  \partial f_1(x) -\{v\}\\
  \partial f_2(x) -\{v\}
\end{pmatrix}
\mbox{ and }
T(x,v) = \begin{pmatrix}
  \partial g_1(x) -\{v\}\\
  \partial g_2(x) -\{v\}
\end{pmatrix}
\]
Clearly, $0\in S(x,v)$ implies that $0 \in \partial f_2(x) - \partial f_1(x)$; and $0 \in \partial f_2(x) - \partial f_1(x)$ implies that there exists a ``multiplier vector'' $v\in\reals^n$ such that $0\in S(x,v)$.
A bound on $\hatsetd_\rho(\gph S, \gph T)$ will then via Thm. \ref{thm:approxgeneralizedequations} furnish a bound on the difference between near-stationary points in the ``primal-dual'' space $\reals^n\times \reals^n$ as one passes from minimizing $f_2-f_1$ to minimizing $g_2-g_1$. For simplicity, we adopt the sup-norm for the remainder of this example. Specifically, we find that for $\rho\in \preals$
\[
\hatsetd_\rho(\gph S, \gph T) \leq  \max_{i=1,2}  \sup_{\|x\|_\infty\leq\rho}  \hatsetd_{2\rho}\big( \partial f_i(x),\partial g_i(x)\big).
\]
To see this let $((\bar x, \bar v), (\bar y_1, \bar y_2)) \in \gph S\cap \ball_{\reals^{4n}}(\rho)$, i.e., $\bar y_1 + \bar v \in \partial f_1(\bar x)$ and  $\bar y_2 + \bar v \in \partial f_2(\bar x)$.
For $i=1, 2$, since $\|\bar x\|_\infty \leq \rho$ and $\|\bar y_i + \bar v \|_\infty \leq 2\rho$, there exists $y_i\in\reals^n$ such that
\[
y_i+\bar v \in \partial g_i(\bar x) \mbox{ and } \|(\bar y_i + \bar v) - (y_i + \bar v)\|_\infty \leq \hatsetd_{2\rho}(\partial f_i(\bar x), \partial g_i(\bar x)),
\]
which implies $((\bar x, \bar v), (y_1, y_2)) \in \gph T$. The distance between $((\bar x, \bar v), (y_1, y_2))$ and $((\bar x, \bar v), (\bar y_1, \bar y_2))$ then yields the stated upper bound on $\hatsetd_\rho(\gph S, \gph T)$.\\

{\it Example 4: KKT conditions}. Theorem \ref{thm:approxgeneralizedequations} also applies to the KKT conditions for the problem
\[
\nnmin f_0(x) \mbox{ subject to } f_i(x) \leq 0 \mbox{ for } i = 1, \dots, m, \mbox{ with smooth } f_i:\reals^n\to \reals,
\]
when compared to those of an alternative, possibly approximating, problem obtained by replacing the functions by the smooth functions $g_0, g_1, \dots, g_m$. Clearly, $(x,y)\in \reals^{n+m}$ satisfies the KKT conditions for the actual problem if and only if $0 \in S(x,y)$ and likewise those of the alternative problem if and only if $0\in T(x,y)$, where the set-valued mappings $S,T:\reals^{n+m} \tto \reals^{3m+n}$ have values
\begin{equation*}\label{eqn:KKTmaps}
S(x,y) = \begin{pmatrix}
[f_1(x), \infty) \\
\vdots\\
[f_m(x), \infty) \\
(-\infty, y_1]\\
\vdots\\
(-\infty, y_m]\\
\{y_1f_1(x)\}\\
\vdots\\
\{y_mf_m(x)\}\\
\{\nabla f_0(x) + \sum_{i=1}^m y_i \nabla f_i(x)\}
\end{pmatrix}
T(x,y) = \begin{pmatrix}
[g_1(x), \infty) \\
\vdots\\
[g_m(x), \infty) \\
(-\infty, y_1]\\
\vdots\\
(-\infty, y_m]\\
\{y_1g_1(x)\}\\
\vdots\\
\{y_mg_m(x)\}\\
\{\nabla g_0(x) + \sum_{i=1}^m y_i \nabla g_i(x)\}
\end{pmatrix}
\end{equation*}
with $y = (y_1, \dots, y_m)$. A bound on the truncated Hausdorff distance between the graphs of these two set-valued mappings furnishes the critical component in the application of Thm. \ref{thm:approxgeneralizedequations}. In this example, we equip $\reals^{n+m}$ and $\reals^{3m+n}$ with the sup-norm.  Then, for $\rho\in \preals$,
\[
\hatsetd_\rho(\gph S, \gph T) \leq \max\{\delta, \rho\delta, (1+m\rho)\eta\},
\]
where
\[
\delta = \max_{i=1, \dots, m}\sup_{\|x\|_\infty \leq \rho} |f_i(x) - g_i(x)| \mbox{ and } \eta = \max_{i=0, \dots, m}\sup_{\|x\|_\infty \leq \rho} \|\nabla f_i(x) - \nabla g_i(x)\|_\infty.
\]
This assertion is realized as follows. Let $((x,y), (u,v,w,s)) \in \gph S \cap \ball_{\reals^{4m+2n}}(\rho)$ be arbitrary and construct $\bar x = x$, $\bar y = y$, $\bar u = (\bar u_1, \dots, \bar u_m)$, with $\bar u_i = \max\{g_i(x), u_i\}$ for all $i$, $\bar v = v$, $\bar w = (\bar w_1, \dots, \bar w_m)$, with $\bar w_i = y_i g_i(x)$ for all $i$, and $\bar s = \nabla g_0(x) + \sum_{i=1}^m y_i \nabla g_i(x)$. It is trivial to verify that $((\bar x,\bar y), (\bar u,\bar v,\bar w,\bar s)) \in \gph T$. For all $i$,
\[
|u_i - \bar u_i| \leq \begin{cases}
0 &\mbox{ if } u_i \geq g_i(x)\\
\delta &\mbox{ otherwise}
\end{cases}
\]
\[
|w_i - \bar w_i| \leq |y_i|\big|f_i(x) - g_i(x)\big| \leq \rho\delta
\]
\[
\|s-\bar s\|_\infty\leq  \|\nabla f_0(x) - \nabla g_0(x)\|_\infty + \sum_{i=1}^m |y_i|\|\nabla f_i(x) - \nabla g_i(x)\|_\infty \leq (1+m\rho)\eta.
\]
Consequently, the distance between $((x,y), (u,v,w,s))$ and $((\bar x,\bar y), (\bar u,\bar v,\bar w,\bar s))$ is at most $\max\{\delta,$ $\rho\delta,$ $(1+m\rho)\eta\}$ and we have that $\exs((\gph S \cap \ball_{\reals^{4m+2n}}(\rho); \gph T))$ is bounded by the same quantity. The assertion then follows by symmetry.

We see that despite the fact that minimizers of inequality-constrained problems are {\it unstable} under pointwise perturbations of the constraint functions (cf.  Section 4), the KKT system has {\it stable} solutions in the sense that the excess of near-solutions of one KKT system over those of the other exhibits a Lipschitz property in those perturbations.\\

We end the paper with a result that generalizes the ideas of Examples 3 and 4.  For a proper lsc function $\phi:\reals^m\to \Reals$ and a smooth mapping $F:\reals^n\to \reals^m$, we recall that under rather weak assumptions\footnote{For example, if $\phi$ is convex, then it suffices that $\dom \phi$ cannot be separated from the range of the linearized mapping $w\mapsto F(\bar x) + \nabla F(\bar x)w$ for a local minimizer $\bar x$.} the composite function $\phi\circ F$ has $0 \in \nabla F(x)^\top \partial \phi(F(x))$ as a necessary optimality condition \cite[Thm. 10.6]{VaAn}, where the $m\times n$-matrix $\nabla F(x)$ is the Jacobian of $F$ at $x$. By introducing auxiliary vectors $y,z\in\reals^m$, the optimality condition is equivalently stated in terms of the set-valued mapping $S:\reals^n\times\reals^m\times\reals^m\tto \reals^m\times\reals^m\times\reals^n$ as $0 \in S(x,y,z)$, with
\begin{equation}\label{eqn:compositionMapping}
  S(x,y,z) = \begin{pmatrix}
    \{F(x) - z\}\\
    \partial \phi(z) - \{y\}\\
    \{\nabla F(x)^\top y\}
  \end{pmatrix}.
\end{equation}
Since $0\in S(x,y,z)$ is also an optimality condition for the problem of minimizing $\phi(z)$ subject to $F(x) = z$, $y$ can be interpreted as a multiplier vector and $z$ as representing feasibility. Parallel conditions hold for a composite function $\psi\circ G$ expressed in terms of $\psi:\reals^m\to \Reals$ and $G:\reals^n\to\reals^m$, which we may think of as approximations of $\phi$ and $F$. Specifically, under the appropriate assumptions, an optimality condition becomes $0\in T(x,y,z)$, where the set-valued mapping $T:\reals^n\times\reals^m\times\reals^m\tto \reals^m\times\reals^m\times\reals^n$ has
\begin{equation}\label{eqn:compositionMapping2}
  T(x,y,z) = \begin{pmatrix}
    \{G(x) - z\}\\
    \partial \psi(z) - \{y\}\\
    \{\nabla G(x)^\top y\}
  \end{pmatrix}.
\end{equation}
In view of Thm. \ref{thm:approxgeneralizedequations}, a bound on $\hatsetd_\rho(\gph S, \gph T)$ leads to estimates of the change in near-stationary points as we pass from $\phi\circ F$ to $\psi\circ G$.

\begin{theorem}{\rm (stationarity of composite functions).} For proper lsc functions $\phi,\psi:\reals^m\to \Reals$, smooth mappings $F,G:\reals^n\to \reals^m$, and the resulting set-valued mappings $S$ and $T$ expressed in  \eqref{eqn:compositionMapping} and \eqref{eqn:compositionMapping2}, we have for
$\rho\in \preals$ that\footnote{Here, $\hatsetd_\rho$ is defined in terms of the product norm on $\reals^n\times\reals^m\times\reals^m\times\reals^m\times\reals^m\times\reals^n$ constructed by any norms on $\reals^n$ and $\reals^m$ and the matrix norm is any one compatible with the norm on $\reals^m$.}
\begin{align*}
\hatsetd_\rho(\gph S, \gph T) \leq \sup_{\|x\|\leq\rho} \max\Big\{ &\big\|G(x) - F(x)\big\| + \hatsetd_{2\rho} \big(\gph \partial \phi, \gph \partial \psi\big),\\
& \rho\big\|\nabla G(x)^\top - \nabla F(x)^\top\big\| \Big\}.
\end{align*}
\end{theorem}
\state Proof. Suppose that $((\bar x, \bar y, \bar z), (\bar u, \bar v, \bar w)) \in \gph S \cap \ball_{X}(\rho)$, where $X = \reals^n\times\reals^m\times\reals^m\times\reals^m\times\reals^m\times\reals^n$ and using the norm indicated in the footnote. Then,
\[
\bar u = F(\bar x) -\bar z, ~ \bar v + \bar y \in \partial \phi(\bar z), ~\bar w = \nabla F(\bar x)^\top \bar y.
\]
Since $(\bar z, \bar v + \bar y) \in \gph \partial \phi \cap \ball_{\reals^{m}\times\reals^m}(2\rho)$ (using the product norm on $\reals^m \times \reals^m$) and the fact that $\gph\partial \psi$ is nonempty \cite[Cor. 8.10]{VaAn}, there exist $z,v\in\reals^m$ such that $(z,v+\bar y) \in \gph \partial \psi$ and neither $\|z-\bar z\|$ nor $\|(\bar v - \bar y) - (v-\bar y)\|$ exceed $\hatsetd_{2\rho}(\gph \partial \phi, \gph \partial \psi)$. Construct $u = G(\bar x) - z$ and $w = \nabla G(\bar x)^\top \bar y$. Clearly, $((\bar x, \bar y, z), (u, v, w)) \in \gph T$ and
\[
\|u-\bar u\| = \big\|(G(\bar x) - z) - (F(\bar x) - \bar z)\big\| \leq \big\|G(\bar x) - F(\bar x)\big\| + \hatsetd_{2\rho}(\gph \partial \phi, \gph \partial \psi).
\]
Moreover, due to the assumed compatibility of the adopted matrix norm relative to the norm on $\reals^m$,
\[
\|w - \bar w\| = \big\|\nabla G(\bar x)^\top\bar y - \nabla F(\bar x)^\top \bar y\big\| \leq \rho\big\|\nabla G(\bar x)^\top - \nabla F(\bar x)^\top\big\|.
\]
The point $((\bar x, \bar y, z), (u, v, w))$ is therefore within a distance of
\[
\max\Big\{\big\|G(\bar x) - F(\bar x)\big\| + \hatsetd_{2\rho}(\gph \partial \phi, \gph \partial \psi)  , ~\rho\big\|\nabla G(\bar x)^\top - \nabla F(\bar x)^\top\big\| \Big\}
\]
of $((\bar x, \bar y, \bar z), (\bar u, \bar v, \bar w))$, which establishes the conclusion after we realize the  obvious symmetry in the result.\eop

\appendix
\section{Proofs}

Proof of Prop. \ref{prop:solestimates}. Denote by $d_X$ the metric on $X$ and $\eta = \hatsetd_\rho(\epi f, \epi g)$. Let $\gamma \in (0, \rho - \epsilon - \inf f)$. Since $\gamma\mbox{-}\nargmin f \cap \ball_X(\rho) \neq \emptyset$, there exists $\bar x\in \ball_X(\rho)$ such that $f(\bar x) \leq \inf f + \gamma < \rho - \epsilon\leq \rho$. Moreover, $f(\bar x) \geq \inf f \geq -\rho$. Thus, $(\bar x, f(\bar x)) \in \epi f \cap \ball_{X\times \reals}(\rho)$ and there exists $(x,\alpha) \in \epi g$ such that $\max\{d_X(x, \bar x), |\alpha-f(\bar x)|\} \leq \dist((\bar x, f(\bar x)), \epi g) + \gamma$. Then,
\[
\eta \geq \exs\big(\epi f \cap \ball_{X\times \reals}(\rho); \epi g\big) \geq \dist\big((\bar x, f(\bar x)), \epi g\big) \geq  d_X\big(x, \bar x\big) - \gamma
\]
and also $\eta \geq |\alpha - f(\bar x)| -\gamma$. Collecting the above results yield $\inf g \leq g(x) \leq \alpha \leq f(\bar x) + \eta + \gamma \leq \inf f + \eta + 2\gamma$. Since $\gamma$ is arbitrary, we have established that $\inf g \leq \inf f + \eta$. The same argument with the roles of $f$ and $g$ reversed leads to the first conclusion.

Let $\bar x \in \epsilon\mbox{-}\nargmin g \cap \ball_X(\rho)$. Then, $g(\bar x) \leq \inf g + \epsilon < \rho$, $g(\bar x) \geq \inf g \geq -\rho$, and $(\bar x, g(\bar x)) \in \epi g \cap \ball_{X\times \reals}(\rho)$. Let $\gamma>0$. There exists $(x,\alpha) \in \epi f$ such that $\max\{d_X(x, \bar x), |\alpha-g(\bar x)|\} \leq \dist((\bar x, g(\bar x)), \epi f) + \gamma$. Consequently, $\eta \geq d_X(x, \bar x) - \gamma \mbox{ and } \eta \geq |\alpha - g(\bar x)| -\gamma$. These facts together with the first conclusion establish that $f(x) \leq \alpha \leq g(\bar x) + \eta + \gamma \leq \inf g + \epsilon + \eta + \gamma \leq \inf f + \epsilon + 2\eta + \gamma$. Thus, $x\in (\epsilon + 2\eta + \gamma)\mbox{-}\nargmin f$ and $d_X(x,\bar x) \leq \eta + \gamma$, and then also $\exs(\epsilon\mbox{-}\nargmin g \cap \ball_X(\rho); ~(\epsilon+2\eta+\bar\gamma)\mbox{-}\nargmin f\big) \leq \eta+ \gamma$ when $\bar \gamma \geq \gamma$. Since $\gamma$ is arbitrary, the second conclusion follows.\eop\\

Proof of Prop. \ref{prop:levelsets}. Let $\bar x \in \nlev_\delta g \cap \ball_X(\rho)$ and $B = \ball_{X\times \reals}(\rho)$.  Then, $g(\bar x) \leq \delta \leq \rho$. There are two cases. Suppose that $g(\bar x) \geq -\rho$. Then, $(\bar x, g(\bar x)) \in \epi g \cap B$. Let $\gamma \in (0, \epsilon - \delta - \exs(\epi g \cap B; \epi f))$. There exists $(x,\alpha) \in \epi f$ such that $\max\{d_X(x, \bar x), |\alpha - g(\bar x)|\} \leq \dist((\bar x, g(\bar x)), \epi f) + \gamma \leq \exs\big(\epi g \cap B; \epi f\big) + \gamma$. Consequently,
\[
f(x) \leq \alpha \leq g(\bar x) + \exs\big(\epi g \cap B; \epi f\big) + \gamma \leq \delta + \exs\big(\epi g \cap B; \epi f\big) + \gamma < \epsilon.
\]
Thus, $x\in \nlev_\epsilon f$ and $d_X(x,\bar x) \leq \exs(\epi g \cap B; \epi f)) + \gamma$. This implies that
\[
\exs(\nlev_\delta g \cap \ball_X(\rho); \nlev_{\epsilon} f) \leq \exs(\epi g \cap B; \epi f) + \gamma.
\]
If $g(\bar x) < - \rho$, the same holds because the arguments in that case can be carried out with $g(\bar x)$ replaced by $-\rho$.
Since $\gamma$ is arbitrary, the second conclusion follows.\eop\\

Proof of Prop. \ref{prop:sumofsets}. Let $C = \sum_{i=1}^m C_i$, $D = \sum_{i=1}^m D_i$, and $\epsilon > 0$. Suppose without loss of generality that $\hatsetd_\rho(C,D) = \exs(C\cap \ball_X(\rho); D)$. If $C\cap \ball_X(\rho)=\emptyset$, $\hatsetd_\rho(C,D) = 0$ and the result holds trivially. Thus, suppose that $C\cap \ball_X(\rho)\neq\emptyset$. Then, there are $x_i\in C_i$ and $y_i\in D_i$, $i=1, \dots, m$, such that $x = \sum_{i=1}^m x_i \in C\cap \ball_X(\rho)$, $\|x_i-y_i\| \leq \dist(x_i,D_i) + \epsilon$, and
\[
\hatsetd_\rho(C,D) \leq \dist(x,D) + \epsilon \leq \|x-y\| + \epsilon \leq \sum_{i=1}^m \|x_i-y_i\| + \epsilon \mbox{ where } y=\sum_{i=1}^m y_i.
\]
Since $x_i\in C_i$ implies $x_i\in\ball_X(\rho)$,
\[
\|x_i - y_i\| \leq \dist(x_i,D_i) + \epsilon \leq \exs(C_i\cap\ball_X(\rho); D_i) + \epsilon\leq \hatsetd_\rho(C_i,D_i)+ \epsilon.
\]
Hence, $\hatsetd_\rho(C,D) \leq \sum_{i=1}^m \hatsetd_\rho(C_i,D_i) + (m+1)\epsilon$. Since $\epsilon$ is arbitrary, the first conclusion follows. Under the relaxed assumption, $x_1\in \ball_X(m\rho)$ because $\{x,x_i\in\ball_X(\rho), i=2, \dots, m\}$.
Thus,
\[
\|x_1 - y_1\| \leq \dist(x_1,D_1) + \epsilon \leq \exs(C_1\cap\ball_X(m\rho); D_1) + \epsilon\leq \hatsetd_{m\rho}(C_1,D_1)+ \epsilon.
\]
Since the other arguments carry over, the second conclusion follows.\eop\\

Proof of Prop. \ref{prop:altform}. Let $\eta = \hatsetd_\rho(\epi f, \epi g)$ and $\epsilon>0$. Suppose that $(x, f(x))\in \epi f \cap \ball_{X\times \reals}(\rho)$. Then, there exist $(\bar x,\bar \alpha) \in \epi g$ such that $d_X(\bar x,x) \leq \eta + \epsilon$, $|\alpha - f(x)|\leq \eta + \epsilon$, and $g(\bar x) \leq \alpha< \infty$. Thus, $g(\bar x) \leq \alpha \leq f(x) + \eta + \epsilon \leq \max\{f(x), -\rho\} + \eta + \epsilon$. This establishes that $\inf_{\ball(x,\eta+\epsilon)} g \leq \max\{f(x), -\rho\} + \eta + \epsilon$ for $x\in \nlev_\rho f \cap \ball_X(\rho)$ and $f(x) \geq -\rho$. Suppose that $x\in \nlev_\rho f \cap \ball_X(\rho)$ and $f(x) < -\rho$. Then, $(x, -\rho)\in \epi f \cap \ball_{X\times \reals}(\rho)$ and there exist $(\bar x,\bar \alpha) \in \epi g$ such that $d_X(\bar x,x) \leq \eta + \epsilon$, $|\alpha + \rho|\leq \eta + \epsilon$, and $g(\bar x) \leq \alpha< \infty$. Consequently,
\[
\ninf_{\ball_X(x,\eta+\epsilon)} g \leq g(x) \leq \alpha \leq -\rho + \eta + \epsilon \leq \max\{f(x), -\rho\} + \eta + \epsilon.
\]
Repeating the arguments with the roles of $f$ and $g$ reversed, we establish that the two sets of constraint on the right-hand side in the proposition is satisfied with $\eta+\epsilon$. Thus, the right-hand side does not exceed $\eta+ \epsilon$. Since $\epsilon$ is arbitrary, the right-hand side furnishes a lower bound on $\hatsetd_\rho(\epi f,\epi g)$. By \cite[Prop. 3.2]{Royset.18}, it is also an upper bound; the lsc assumption in that proposition is not needed in its proof.\eop

\noindent {\bf Acknowledgement.} This work is supported in part by DARPA (Lagrange) under HR0011-8-34187, ONR (Science of Autonomy) under  N0001419WX00183, and AFOSR (Optimization and Discrete Mathematics) under  F4FGA08272G001.

\bibliographystyle{plain}
\bibliography{refs}

\end{document}